\documentclass[final,3p,12pt, times]{elsarticle}

\usepackage[utf8]{inputenc}
\usepackage{amsmath,amssymb,amsthm,mathtools}
\usepackage{mathptmx}
\usepackage{courier}
\usepackage{microtype}
\usepackage{enumitem}
\usepackage{needspace}
\usepackage{geometry}
\usepackage{titlesec}
\usepackage{hyperref}
\usepackage[nameinlink,noabbrev]{cleveref}

\hypersetup{
  hidelinks,
  pdftitle={Hausdorff reflections and obstruction ideals for filter topologies on MV-algebras},
  pdfauthor={Yanyan Dong, Bin Pang, Yi Shi}
}
\biboptions{sort&compress}

\setlist{leftmargin=2em}
\setlist[enumerate,1]{label=(\arabic*),leftmargin=0pt,labelindent=0pt,labelwidth=1.4em,labelsep=0.4em,itemindent=1.8em,align=left,font=\upshape,topsep=2pt,itemsep=1pt,parsep=0pt,partopsep=0pt}
\setlist[itemize,1]{leftmargin=2em,topsep=2pt,itemsep=1pt,parsep=0pt,partopsep=0pt}
\titleformat{\section}{\normalfont\bfseries\large}{\thesection.}{0.55em}{}
\titleformat{\subsection}{\normalfont\bfseries\normalsize}{\thesubsection.}{0.55em}{}
\titlespacing*{\section}{0pt}{2.4ex plus .6ex minus .2ex}{1.25ex plus .2ex}
\titlespacing*{\subsection}{0pt}{2.0ex plus .5ex minus .2ex}{1.0ex plus .2ex}

\newtheoremstyle{jpaaplain}
  {6pt}{6pt}{\itshape}{}{\bfseries}{.}{0.5em}{}
\newtheoremstyle{jpaadefinition}
  {6pt}{6pt}{\normalfont}{}{\bfseries}{.}{0.5em}{}
\newtheoremstyle{jpaaremark}
  {6pt}{6pt}{\normalfont}{}{\bfseries}{.}{0.5em}{}

\theoremstyle{jpaaplain}
\newtheorem{theorem}{Theorem}[section]
\newtheorem*{maintheorem}{Main theorem}
\newtheorem{lemma}[theorem]{Lemma}
\newtheorem{proposition}[theorem]{Proposition}
\newtheorem{corollary}[theorem]{Corollary}

\theoremstyle{jpaadefinition}
\newtheorem{definition}[theorem]{Definition}
\newtheorem{example}[theorem]{Example}

\theoremstyle{jpaaremark}
\newtheorem{remark}[theorem]{Remark}

\newcommand{\Cint}{\mathcal C_{\mathrm{int}}}
\newcommand{\Cd}{\mathcal C_d}

\newcommand{\co}{\operatorname{co}}

\makeatletter
\def\ps@pprintTitle{\ps@plain}
\makeatother
\usepackage{indentfirst}
\newif\ifincludeconvexity
\includeconvexityfalse 

\begin{document}
\begin{frontmatter}

\title{Hausdorff reflections and obstruction ideals for filter topologies on MV-algebras}

\author[1]{Yanyan Dong}
\ead{dongyanyan187@163.com}
\author[2]{Bin Pang\corref{cor1}}
\ead{pangbin1205@163.com}
\cortext[cor1]{Corresponding author.}

\author[3]{Yi Shi}
\ead{shiyi7353@126.com}

\address[1]{School of Mathematics and Statistics, Yancheng Teachers University, Yancheng, Jiangsu 224002, PR China}
\address[2]{Beijing Key Laboratory on MCAACI, School of Mathematics and Statistics, Beijing Institute of Technology, Beijing 102488, PR China}

\address[3]{School of Mathematics and Physics, School of Cryptography,
Nanjing Institute of Technology, Nanjing, Jiangsu 211167, PR
China}

\begin{abstract}
MV-algebras arise in many-valued logic and are closely related to
lattice-ordered groups. We characterize Hausdorff separation in
filter topologies through their associated ideals. For every proper
lattice filter, $T_0$, $T_1$ and $T_2$ are equivalent, and the
zero-neighborhood ideal determines a canonical Hausdorff reflection.
For a radius filter satisfying a refinement condition, the filter
topology agrees with the ambient subspace topology exactly when every
radius is a weak order unit. This answers a question of Luan, Zhao
and Yang without requiring divisibility of the representing group.
These results reveal two successive algebraic obstructions: a polar obstruction controlling compatibility with the ambient topology, followed by a relative infinitesimal obstruction controlling the residual failure of Hausdorff separation.
 For a
unital lattice-ordered subring of $C(X,\mathbb R)$, where $X$ is
compact Hausdorff and the unit is $\mathbf1$, the associated exact
sequence admits a positive additive section exactly when the relative
infinitesimal ideal vanishes.
\end{abstract}

\end{frontmatter}

\section{Introduction}
MV-algebras provide an algebraic setting for many-valued logic.
The Mundici equivalence represents them as unit intervals in abelian
lattice-ordered groups with a strong unit
\cite{Chang1958,Chang1959,Mundici1986,Cignoli2000,Mundici2011}.
Their topological properties were studied by Hoo and Weber
\cite{Hoo1997,Weber2012}. In connection with the $C$-topologies of
lattice-ordered groups \cite{Gusic1998,Yang2009}, Luan and Yang
introduced filter topologies on MV-algebras \cite{LuanYang2017}.
For a proper lattice filter $F$ of an MV-algebra $A$, these
topologies have basic neighborhoods
\[
 B_F(x,r)=\{y\in A:r\ominus d(x,y)\in F\},\qquad r\in F,
\]
where $d$ is Chang's distance. Further properties of these topologies
were studied in \cite{Asadzadeh2020,WuLuanYang2020,LuanWeberYang2021}.

Our first question is whether $T_0$ already implies Hausdorffness
for every filter topology. Luan, Weber and Yang showed that such
a topology makes an MV-algebra quasitopological, but addition need
not be jointly continuous \cite[Corollary~3.2 and Example~3.4]{LuanWeberYang2021}.
For quasitopological MV-algebras, Najafi proved that $T_0$ and
$T_1$ are equivalent, while his equivalence of $T_0$ and $T_2$
requires joint continuity \cite[Theorems~3.3 and 3.5]{Najafi2023}.
The equivalence of $T_0$ and $T_2$ for topological MV-algebras
also appears in Hoo's earlier work \cite{Hoo1997}.
See also \cite{Najafi2017}. Thus the result for topological
MV-algebras does not directly settle the filter-topology case.
The question was recorded as open by Gan, Luan, Deng and Yang
\cite{GanLuanDengYang2024}. We give an affirmative answer in
Theorem~\ref{thm:T0T2}, without assuming joint continuity.

To study separation in a filter topology, we consider the
intersection of the zero-neighborhoods
\[
 N_F=\bigcap_{r\in F}B_F(0,r)
     =\{a\in A:r\ominus a\in F\text{ for every }r\in F\}.
\]
It is an MV-ideal and $d(x,y)\in N_F$ exactly when $x$ and $y$
are topologically indistinguishable. The ideal itself is already
implicit in Bailey's kernel $K(F)$, since $N_F=\{z^*:z\in K(F)\}$
\cite[Theorems~2.1 and 2.3]{Bailey2009}. Bailey's saturation theorem
also identifies the ideals whose quotients preserve membership in
$F$ \cite[Theorem~2.6]{Bailey2009}. The zero-neighborhood ideal
appears in the chain analysis of Luan, Weber and Yang
\cite[Lemma~4.1]{LuanWeberYang2021}.

Let $q:A\longrightarrow A/N_F$ be the quotient map. In
Theorem~\ref{thm:reflection}, we show that $\overline F=q(F)$ is
a proper lattice filter and that its balls pull back to the
corresponding $F$-balls. The quotient $A/N_F$ is Hausdorff, and
every continuous MV-homomorphism from $A$ to a Hausdorff topologized
MV-algebra factors uniquely through $q$. Thus $q$ is the canonical
Hausdorff reflection. We also determine the intermediate quotients
that preserve membership in $F$.

The representation $A=\Gamma(G,u)$ leads to a second question.
Luan, Zhao and Yang observed that a filter topology on a unit
interval can differ from the subspace topology induced by the
$C$-topology on $G$. They asked when the two agree
\cite[Remark~3.1 and Question~2]{LuanZhaoYang2021}.
We consider a nonempty proper lattice filter $C\subseteq G^+$ with
$F=C\cap[0,u]$ and the following radius refinement property
\[
 \text{for every }c\in C\text{ there is }d\in C\text{ with }2d\le c.
\]
Under this condition, Theorem~\ref{thm:weak-unit} answers their
question by showing that the two topologies agree if and only if
every element of $C$ is a weak order unit. Equivalently, the
directed union of principal polars
\[
 D_C=\bigcup_{c\in C}c^\perp
\]
is zero. This includes the usual halving-closed admissible filters
in 2-divisible groups, but requires no divisibility of $G$.

Let $\mathcal N_C$ be the $\ell$-ideal corresponding to $N_F$.
Then $D_C\subseteq\mathcal N_C$ and
$N_F=\mathcal N_C\cap[0,u]$. In Theorem~\ref{thm:decomposition},
we show that $\mathcal N_C/D_C$ is the ideal of infinitesimals
relative to the image filter on $G/D_C$. This ideal determines
whether the quotient topology is Hausdorff. At the MV-algebra level, the corresponding polar obstruction ideal
$D_F$ is the least kernel whose quotient makes all image radii weak
order units. The two obstruction ideals are stable under every
intermediate quotient that preserves the filter, yielding explicit
quotient formulas and criteria for both subspace compatibility and
Hausdorffness.

The ambient group topology belongs to the positive-filter theory of
Kopperman, Pajoohesh and Richmond \cite{KoppermanPajooheshRichmond2011},
with completions studied in \cite{JordanPajoohesh2018}. Related work
on topological MV-algebras includes initial topologies and
isomorphism theorems \cite{YinXieYang2025}, uniformization and
quotient homomorphisms \cite{LiYang2026,XieYang2026,XieYangUniform2026},
and topological Mundici representations
\cite{JiaZhaoYang2026,MengShenZhang2026}. The last of these also
treats a quasitopological correspondence.

The following theorem summarizes the separation, reflection and
structural results proved below.

\begin{maintheorem}
Let $F$ be a proper lattice filter of an MV-algebra $A$, and let
$N_F$ be the intersection of its zero-neighborhoods.
\begin{enumerate}
\item The separation axioms $T_0$, $T_1$ and $T_2$ are equivalent
for $\tau_F$. The image $\overline F$ of $F$ under the quotient map
$q:A\longrightarrow A/N_F$ is a proper lattice filter and
\[
 q^{-1}\bigl(B_{\overline F}(q(x),q(r))\bigr)=B_F(x,r)
 \qquad(x\in A,\ r\in F).
\]
The quotient $A/N_F$, equipped with $\tau_{\overline F}$, is
Hausdorff and $q$ is the canonical Hausdorff reflection for
continuous MV-homomorphisms.
\item Suppose that $A=\Gamma(G,u)$ and $F=C\cap[0,u]$, where $C$
is a nonempty proper lattice filter in $G^+$ satisfying the radius
refinement condition above. Put
\[
 D_C=\bigcup_{c\in C}c^\perp,\qquad
 \mathcal N_C=\{g\in G:(c-|g|)^+\in C\text{ for every }c\in C\}.
\]
Write $\tau_C^G$ for the ambient $C$-topology on $G$. Then
$N_F=\mathcal N_C\cap[0,u]$, and
\[
 \tau_F=\tau_C^G|_{[0,u]}
 \quad\Longleftrightarrow\quad D_C=\{0\}
 \quad\Longleftrightarrow\quad
 \text{every }c\in C\text{ is a weak order unit}.
\]
For $p:G\longrightarrow H=G/D_C$ and $E=p(C)$, there is an
exact sequence of $\ell$-groups
\[
 0\longrightarrow D_C\longrightarrow\mathcal N_C
 \xrightarrow{\ p\ }\operatorname{Inf}_E(H)\longrightarrow0,
\]
where $\operatorname{Inf}_E(H)=\{h\in H:n|h|\le e\text{ for all }
e\in E,\ n\ge1\}$.
\end{enumerate}
\end{maintheorem}

The proof of the separation equivalence uses intersections of
filter balls. The polar quotient relates these balls to the ambient
group topology. We also determine the obstruction ideals for finite
products of chains and for groups of continuous functions. For
function rings, pointwise squares are used to study positive sections
of the exact sequence. For the vanishing ideals in these examples,
see \cite{GillmanJerison1976,HenriksenLarsonMartinezWoods1994}.

Section~\ref{sec:prelim} recalls some basic notions and results.
Sections~\ref{sec:separation} and \ref{sec:structure} establish the
reflection and the quotient laws in MV-algebras and $\ell$-groups,
respectively. Section~\ref{sec:models} proves the product
classification and the function-group results.

\section{Preliminaries}\label{sec:prelim}

We recall the notation and results needed below. Unless otherwise stated, $A$ denotes a nontrivial MV-algebra, $G$ an abelian $\ell$-group, and $n$ a positive integer. All lattice filters are assumed to be proper. We follow \cite{AndersonFeil1988,Cignoli2000,Mundici2011}.

\Needspace{10\baselineskip}
\begin{definition}[\cite{Chang1958,Cignoli2000,Mundici2011}]\label{def:mv}
An \emph{MV-algebra} is an algebra $(A,\oplus,{}^*,0)$ satisfying the following conditions for all $x,y\in A$.
\begin{enumerate}
\item $(A, \oplus, 0)$ is a commutative monoid.
\item $\left(x^*\right)^*=x$.
\item $x \oplus 0^*=0^*$.
\item $\left(x^* \oplus y\right)^* \oplus y=\left(y^* \oplus x\right)^* \oplus x$.
\end{enumerate}
Write $1=0^*$ and define
\[
 x\odot y=(x^*\oplus y^*)^*,\qquad
 x\ominus y=x\odot y^*,\qquad
 x\to y=x^*\oplus y.
\]
The natural order is defined by $x\le y$ if and only if $x\to y=1$. It makes $A$ a bounded distributive lattice, where
\[
 x\vee y=(x\ominus y)\oplus y,\qquad
 x\wedge y=(x^*\vee y^*)^*.
\]
\end{definition}

\begin{definition}[\cite{Cignoli2000,Mundici2011}]\label{def:mv-chain}
An MV-algebra whose natural order is total is called an \emph{MV-chain}. The standard MV-chain is $[0,1]$, with
\[
 x\oplus y=\min\{x+y,1\},\qquad x^*=1-x.
\]
Its finite subalgebras are the chains $\text{\L}_{m+1}=\{0,1/m,\ldots,1\}$, $m\ge1$.
\end{definition}

We next recall the representation of MV-algebras by lattice-ordered groups.

\begin{definition}[\cite{AndersonFeil1988}]\label{def:lgroup}
An abelian \emph{lattice-ordered group}, or \emph{$\ell$-group}, is an abelian group $G$ endowed with a translation-invariant lattice order. We write
\[
 G^+=\{g\in G:g\ge0\},\qquad
 g^+=g\vee0,\qquad g^-=(-g)\vee0,\qquad |g|=g\vee(-g).
\]
A positive element $u$ is a \emph{strong unit} if for every $g\in G$ there exists $n\ge1$ such that $|g|\le nu$.
\end{definition}

The representation theorem relates MV-algebras to lattice-ordered
groups \cite{Chang1959,Mundici1986,Cignoli2000}.

\begin{theorem}\label{thm:representation-prelim}
Every MV-algebra embeds subdirectly into a product of MV-chains and is isomorphic to
$
 \Gamma(G,u)=[0,u]
$
for an abelian $\ell$-group $G$ with strong unit $u$. In this representation,
\[
 x\oplus y=(x+y)\wedge u,\qquad
 x^*=u-x,\qquad
 x\ominus y=(x-y)^+,
\]
and the lattice operations are those of $G$. For an MV-chain, the representing $\ell$-group is linearly ordered.
\end{theorem}

By Theorem~\ref{thm:representation-prelim}, identities and finite systems
of order inequalities can be checked in the chain factors.

\begin{definition}[\cite{Cignoli2000,Mundici2011}]\label{def:chang-distance}
For an MV-algebra $A$, the \emph{Chang distance} is the map $d:A\times A\longrightarrow A$ defined by
\[d(x,y)=(x\ominus y)\oplus(y\ominus x)\qquad(x,y\in A).\]
\end{definition}

The following properties of Chang's distance are standard
\cite{Cignoli2000,Mundici2011}.

\begin{proposition}\label{prop:lattice-distance}
Let $A$ be an MV-algebra. Then
\[d(x,y)=(x\ominus y)\vee(y\ominus x)\]
for all $x,y\in A$.
\end{proposition}

\Needspace{8\baselineskip}
\begin{proposition}\label{prop:distance-identities-prelim}
Let $A$ be an MV-algebra and $x,y,z\in A$. Then the following statements hold.
\begin{enumerate}
\item $d(x,y)=d(y,x)$.
\item $d(x,y)=0$ if and only if $x=y$.
\item $d(x,z)\le d(x,y)\oplus d(y,z)$.
\item $x\ominus(y\vee z)=(x\ominus y)\wedge(x\ominus z)$.
\item $(x\ominus y)\ominus z=x\ominus(y\oplus z)$.
\item If $A=\Gamma(G,u)$, then $d(x,y)=|x-y|$.
\end{enumerate}
\end{proposition}

Chang's distance also describes the congruence associated with an
MV-ideal.

\begin{definition}[\cite{Cignoli2000,Mundici2011}]\label{def:mv-ideal}
A nonempty lower set $I\subseteq A$ closed under $\oplus$ is an \emph{MV-ideal}. The associated congruence is
\[x\equiv_I y\quad\Longleftrightarrow\quad d(x,y)\in I.\]
We write $q_I:A\longrightarrow A/I$ for the quotient homomorphism.
\end{definition}

We also use ideals in the representing $\ell$-group.

\begin{definition}[\cite{AndersonFeil1988}]\label{def:lideal}
An \emph{$\ell$-ideal} $J$ of an abelian $\ell$-group is a solid subgroup, meaning that $h\in J$ and $|g|\le |h|$ imply $g\in J$. Equivalently, $J$ is an order-convex subgroup closed under lattice operations. For the quotient map,
\[
 g+J\le h+J\quad\Longleftrightarrow\quad (g-h)^+\in J.
\]
\end{definition}

The ideals in the two settings are related by the following
correspondence \cite{Cignoli2000,Mundici2011}.

\begin{proposition}\label{prop:ideal-correspondence}
Let $A=\Gamma(G,u)$, where $G$ is an abelian $\ell$-group with strong
unit $u$. Then the following statements hold.
\begin{enumerate}
\item $J\longmapsto J\cap[0,u]$ is an order isomorphism from the
$\ell$-ideals of $G$ to the MV-ideals of $A$.
\item For every $\ell$-ideal $J$ of $G$,
$A/(J\cap[0,u])\cong\Gamma(G/J,u+J)$.
\end{enumerate}
\end{proposition}

We recall disjointness and weak order units in an $\ell$-group.

\begin{definition}[\cite{AndersonFeil1988}]\label{def:order-units}
An $\ell$-group $G$ is \emph{2-divisible} if every $g\in G$ has a half $g/2$. It is \emph{Archimedean} if $g,h\in G^+$ and $ng\le h$ for every $n\ge1$ imply $g=0$. Elements $g,h\in G$ are \emph{disjoint}, written $g\perp h$, if
\[
 |g|\wedge|h|=0.
\]
For $c\in G^+$, its principal polar is
\[
 c^\perp=\{g\in G:|g|\wedge c=0\}.
\]
A positive element $c$ is a \emph{weak order unit} if $c^\perp=\{0\}$.
\end{definition}

The radical $\operatorname{Rad}(A)$ is the intersection of the maximal
MV-ideals of $A$, and $A$ is \emph{semisimple} if
$\operatorname{Rad}(A)=\{0\}$. The radical is related to the
infinitesimals of the representing group as follows
\cite{Cignoli2000,Mundici2011}.

\begin{proposition}\label{prop:radical-infinitesimal}
Let $A=\Gamma(G,u)$, where $G$ is an abelian $\ell$-group with strong
unit $u$, and put
$\operatorname{Inf}_u(G)=\{g\in G:n|g|\le u\text{ for every }n\ge1\}$.
Then the following hold.
\begin{enumerate}
\item $\operatorname{Rad}(A)=\operatorname{Inf}_u(G)\cap[0,u]$.
\item $A$ is semisimple if and only if $G$ is Archimedean.
\end{enumerate}
\end{proposition}

\Needspace{4\baselineskip}
We now recall the definition of a filter topology.

\begin{definition}[\cite{LuanYang2017,LuanWeberYang2021}]\label{def:filter-topology}
A \emph{lattice filter} of an MV-algebra $A$ is a nonempty upward-closed subset $F\subseteq A$ closed under finite meets. It is \emph{proper} if $0\notin F$. This is weaker than an MV-filter, which is also closed under $\odot$. For $x\in A$ and $r\in F$, define
\[B_F(x,r)=\{y\in A:r\ominus d(x,y)\in F\}.\]
The family $\{B_F(x,r):r\in F\}$ is an open neighborhood base at $x$ for the \emph{filter topology} $\tau_F$.
\end{definition}

\begin{remark}\label{rem:radius-monotonicity}
If $r\le s$ with $r,s\in F$, then
\[
 B_F(x,r)\subseteq B_F(x,s).
\]
\end{remark}

We use the following terminology for continuity of the MV-operations.

\begin{definition}[\cite{Najafi2017,LuanWeberYang2021}]\label{def:topological-mv}
A \emph{quasitopological MV-algebra} has continuous negation and separately continuous MV-addition. A \emph{topological MV-algebra} has continuous negation and jointly continuous MV-addition. Every filter topology is quasitopological. We use the standard separation axioms $T_0,T_1,T_2$. Two points are \emph{topologically indistinguishable} if they belong to exactly the same open sets.
\end{definition}

\section{Hausdorff separation and the canonical reflection}\label{sec:separation}
We begin with the separation question recorded in \cite{GanLuanDengYang2024}. Our first goal is to show that, for an arbitrary proper lattice filter, the separation axioms $T_0$, $T_1$ and $T_2$ are equivalent. We then identify the ideal determined by the zero-neighborhoods and show that quotienting by this ideal produces the canonical Hausdorff reflection of the filter topology. The proof of the separation result uses the following elementary inequality.

\begin{lemma}\label{lem:key}
If $a,b,c\in A$ satisfy $a\le b\oplus b$ and $a\le c\oplus c$, then $a\le b\oplus c$.
\end{lemma}
\begin{proof}
By the subdirect representation theorem, it suffices to work in an
MV-chain. If $b\le c$, then
\[
 a\le b\oplus b\le b\oplus c.
\]
If $c\le b$, the same conclusion follows from
$a\le c\oplus c\le b\oplus c$.
\end{proof}

We use this inequality to prove the following separation result.

\begin{theorem}\label{thm:T0T2}
For every proper lattice filter $F$ of an MV-algebra $A$, the separation axioms $T_0$, $T_1$ and $T_2$ are equivalent for $\tau_F$.
\end{theorem}
\begin{proof}
The implications $T_2\Rightarrow T_1\Rightarrow T_0$ hold for
every topological space, so it is enough to prove
$T_0\Rightarrow T_2$. Suppose that $\tau_F$ is $T_0$ and that
distinct points $x,y$ cannot be separated by disjoint
neighborhoods. Put $\delta=d(x,y)$. There is an open set
containing one of $x,y$ and not the other. By the neighborhood-base
property and the symmetry of $d$, we can choose $r_0\in F$ with
\[
 r_0\ominus\delta\notin F.
\]
By assumption, the balls with centers $x$ and $y$ always meet. Choose
\[
 z_1\in B_F(x,r_0)\cap B_F(y,r_0),\qquad
 w_1=d(x,z_1)\vee d(y,z_1).
\]
Using Proposition~\ref{prop:distance-identities-prelim}, we get
\[
 r_1:=r_0\ominus w_1
 =(r_0\ominus d(x,z_1))\wedge(r_0\ominus d(y,z_1))\in F.
\]
Choose $z_2\in B_F(x,r_1)\cap B_F(y,r_1)$ and put
$w_2=d(x,z_2)\vee d(y,z_2)$. The same meet calculation gives
$r_1\ominus w_2\in F$. The subtraction identity in
Proposition~\ref{prop:distance-identities-prelim} therefore yields
\[
 r_0\ominus(w_1\oplus w_2)=r_1\ominus w_2\in F.
\]
For $i=1,2$, the triangle inequality gives
\[
 \delta\le d(x,z_i)\oplus d(z_i,y)\le w_i\oplus w_i.
\]
By Lemma~\ref{lem:key}, $\delta\le w_1\oplus w_2$. Hence
\[
 r_0\ominus\delta\ge r_0\ominus(w_1\oplus w_2)\in F,
\]
which is a contradiction. Therefore $\tau_F$ is Hausdorff.
\end{proof}

Define the \emph{separation ideal associated with $F$} by
\begin{equation}\label{eq:NF}
 N_F=\bigcap_{r\in F}B_F(0,r)
     =\{a\in A:r\ominus a\in F\text{ for every }r\in F\}.
\end{equation}

The ideal property is already contained in Bailey's kernel
construction. In his notation,
\[
 K(F)=\{z\in A:r\odot z\in F\text{ for every }r\in F\},
 \qquad N_F=\{z^*:z\in K(F)\}.
\]
Since $K(F)$ is an implication filter,
\cite[Theorems~2.1 and 2.3]{Bailey2009} imply that $N_F$ is an
MV-ideal. For the zero-neighborhood ideal in an MV-chain, see also
\cite[Lemma~4.1(a)]{LuanWeberYang2021}.

The topological interpretation follows from the neighborhood base.
Indeed, $d(x,y)\in N_F$ means that $y$ lies in every basic
neighborhood of $x$. Symmetry of $d$ then says that $x$ lies in
every basic neighborhood of $y$ as well. Thus
\[
 x,y\text{ are topologically indistinguishable}
 \quad\Longleftrightarrow\quad d(x,y)\in N_F.
\]
In particular, $\tau_F$ is $T_0$ exactly when $N_F=\{0\}$,
as also follows from \cite[Theorem~3.3]{Najafi2023} by negation.
The next theorem identifies the filter topology on the quotient
and establishes its universal property.

\begin{theorem}\label{thm:reflection}
Let $F$ be a proper lattice filter of $A$, let $N_F$ be the ideal
in \eqref{eq:NF}, and let $q:A\longrightarrow A/N_F$ be the quotient
homomorphism. Put $\overline F=q(F)$. Then the following hold.
\begin{enumerate}
\item The set $\overline F$ is a proper lattice filter,
$q^{-1}(\overline F)=F$, and
\begin{equation}\label{eq:reflection-balls}
 q^{-1}\bigl(B_{\overline F}(q(x),q(r))\bigr)=B_F(x,r)
 \qquad(x\in A,\ r\in F).
\end{equation}
\item The map $q$ is continuous and open. The topology
$\tau_{\overline F}$ is the quotient topology, it is Hausdorff, and
every $\tau_F$-open set is the inverse image of a unique
$\tau_{\overline F}$-open set.
\item Every continuous MV-homomorphism
$f:(A,\tau_F)\longrightarrow B$ into an MV-algebra endowed with a
Hausdorff topology factors uniquely as
$f=\overline f\circ q$, where
$\overline f:(A/N_F,\tau_{\overline F})\longrightarrow B$ is a
continuous MV-homomorphism.
\end{enumerate}
Thus $q$ is the Hausdorff reflection in the category of
filter-topological MV-algebras and continuous MV-homomorphisms.
\end{theorem}
\begin{proof}
The image $\overline F=q(F)$ is closed under finite meets. If
$q(f)\le q(a)$ with $f\in F$, then $q(f\vee a)=q(a)$ and
$f\vee a\in F$, so $\overline F$ is upward closed. Moreover,
$N_F\cap F=\varnothing$, since $a\in N_F\cap F$ would give
$a\ominus a=0\in F$. Thus $0\notin\overline F$, and $\overline F$
is proper.

Bailey's saturation theorem \cite[Theorem~2.6]{Bailey2009},
applied to the implication filter $K(F)=\{a^*:a\in N_F\}$,
gives $q^{-1}(\overline F)=F$. The quotient by this implication
filter is the same quotient as the one by the MV-ideal $N_F$.

Since $q$ preserves Chang's distance and $\ominus$, for $x,y\in A$ and
$r\in F$ we have
\[
\begin{aligned}
 q(y)\in B_{\overline F}(q(x),q(r))
 &\Longleftrightarrow q(r\ominus d(x,y))\in\overline F\\
 &\Longleftrightarrow r\ominus d(x,y)\in F\\
 &\Longleftrightarrow y\in B_F(x,r).
\end{aligned}
\]
Therefore the basic balls satisfy the asserted inverse-image identity.
Surjectivity of $q$ also gives
$q(B_F(x,r))=B_{\overline F}(q(x),q(r))$. Hence $q$ is continuous and
open, and $\tau_{\overline F}$ is the quotient topology.
Every open set in $A$ is a union of basic balls and hence is
$q$-saturated by \eqref{eq:reflection-balls}. It is therefore the
inverse image of its open image under $q$. Uniqueness follows from
surjectivity.

Bailey's quotient formula
\cite[Lemma~5.4]{Bailey2009}, together with
$q^{-1}(\overline F)=F$, gives
$q^{-1}(N_{\overline F})=N_F$ in ideal notation. Since $q$ is
onto, $N_{\overline F}=\{0\}$. Theorem~\ref{thm:T0T2} and the
separation criterion above show that $\tau_{\overline F}$ is
Hausdorff.

If $f:A\longrightarrow B$ is a continuous MV-homomorphism and $B$ is
Hausdorff, then $f(a)=f(0)$ for every $a\in N_F$, since $a$ and $0$
are indistinguishable. Thus $N_F\subseteq\ker f$, and there is a unique
MV-homomorphism $\overline f:A/N_F\longrightarrow B$ with
$f=\overline f\circ q$. For every open set $V\subseteq B$,
\[
 q^{-1}(\overline f^{-1}(V))=f^{-1}(V)
\]
is open in $A$. The quotient property of $q$ implies that
$\overline f^{-1}(V)$ is open, so $\overline f$ is continuous.
This proves the asserted universal property.
\end{proof}

The next result relates intermediate quotients. Its equivalence of
(1) and (2) is Bailey's saturation theorem
\cite[Theorem~2.6]{Bailey2009} in ideal notation. The formula for
the zero-neighborhood ideal follows from \cite[Lemma~5.4]{Bailey2009}.

\begin{theorem}\label{thm:intermediate-quotients}
Let $F$ be a proper lattice filter of an MV-algebra $A$, and let $I$ be an MV-ideal. Write $q_I:A\to A/I$. The following conditions are equivalent.
\begin{enumerate}
\item $I\subseteq N_F$.
\item $q_I^{-1}(q_I(F))=F$.
\end{enumerate}
When these conditions hold, $F_I=q_I(F)$ is a proper lattice filter,
\[
 q_I^{-1}\bigl(B_{F_I}(q_I(x),q_I(r))\bigr)=B_F(x,r),
 \qquad N_{F_I}=q_I(N_F)=N_F/I.
\]
In particular, $q_I$ is continuous and open, every open subset of $A$ is $q_I$-saturated, and the Hausdorff reflection of $A/I$ is canonically $A/N_F$ with its filter topology. Thus the interval $[0,N_F]$ in the lattice of MV-ideals parametrizes exactly the quotients that preserve membership in $F$.
\end{theorem}
\begin{proof}
Put $P=\{a^*:a\in I\}$. This is the implication filter associated
with $I$, and $I\subseteq N_F$ is equivalent to $P\subseteq K(F)$.
The equivalence of (1) and (2) is therefore precisely
\cite[Theorem~2.6]{Bailey2009}. Under these conditions,
\cite[Lemma~5.4]{Bailey2009} gives
\[
 q_I^{-1}(N_{F_I})=N_F,
 \qquad N_{F_I}=q_I(N_F)=N_F/I.
\]

The image $F_I$ is closed under finite meets. If
$q_I(f)\le q_I(a)$ with $f\in F$, then $f\vee a\in F$ and
$q_I(f\vee a)=q_I(a)$, proving upward closure. It is proper because
$I\cap F\subseteq N_F\cap F=\varnothing$.
For $x,y\in A$ and $r\in F$, condition (2) gives
\[
\begin{aligned}
 q_I(y)\in B_{F_I}(q_I(x),q_I(r))
 &\Longleftrightarrow q_I(r\ominus d(x,y))\in F_I\\
 &\Longleftrightarrow r\ominus d(x,y)\in F\\
 &\Longleftrightarrow y\in B_F(x,r).
\end{aligned}
\]
This proves the inverse-image identity. Since $q_I$ is surjective,
it also gives
\[
 q_I(B_F(x,r))=B_{F_I}(q_I(x),q_I(r)).
\]
Every quotient radius belongs to $q_I(F)$, so these identities
prove continuity and openness. They also show that every basic
ball, and hence every open subset of $A$, is $q_I$-saturated.

By the third isomorphism theorem, the canonical map
\[
 \theta:(A/I)/(N_F/I)\longrightarrow A/N_F,
 \qquad \theta((a+I)+(N_F/I))=a+N_F,
\]
is an MV-isomorphism. It takes the image of $F_I$ onto the image
of $F$ and therefore takes basic filter balls onto basic filter
balls. Thus $\theta$ is a homeomorphism. Applying
Theorem~\ref{thm:reflection} to $A/I$ identifies its Hausdorff
reflection with $A/N_F$ as claimed.
\end{proof}

For $r\in A$, write $r^\perp=\{a\in A:a\wedge r=0\}$, and call $r$ a \emph{weak order unit} if $r^\perp=\{0\}$. Define
\[
 D_F=\bigcup_{r\in F}r^\perp.
\]
We call $D_F$ the polar obstruction ideal associated with $F$.
The next theorem describes its universal quotient property.

\begin{theorem}\label{thm:mv-polar-quotient}
Let $F$ be a proper lattice filter of an MV-algebra $A$.
\begin{enumerate}
\item The set $D_F$ is an MV-ideal contained in $N_F$. For $p:A\to A/D_F$, the filter $p(F)$ is proper, every element of $p(F)$ is a weak order unit, and $p$ preserves basic balls by exact inverse image.
\item If $f:A\to B$ is an MV-homomorphism for which every $f(r)$, $r\in F$, is a weak order unit in $B$, then $D_F\subseteq\ker f$. Thus $f$ factors uniquely through $p$. In particular, $D_F$ is the least quotient kernel that makes all image radii weak order units.
\item If $I\subseteq N_F$ and $F_I=q_I(F)$, then
\[
 D_{F_I}=(D_F\vee I)/I,
 \qquad N_{F_I}=N_F/I,
\]
where $\vee$ denotes the join of MV-ideals. Consequently, among the quotients that preserve membership in $F$, those for which all image radii are weak order units are exactly the quotients with kernels in $[D_F,N_F]$. The only Hausdorff quotient in this class has kernel $N_F$.
\end{enumerate}
\end{theorem}
\begin{proof}
Principal polars are $\ell$-ideals in the representing group
\cite{AndersonFeil1988}, so their intersections with $A$ are
MV-ideals by Proposition~\ref{prop:ideal-correspondence}.
The ideals $r^\perp$, $r\in F$, form a directed family because
\[
 r^\perp\cup s^\perp\subseteq(r\wedge s)^\perp,
 \qquad r\wedge s\in F.
\]
Their union $D_F$ is consequently an MV-ideal.
If $a\wedge r=0$ for some $r\in F$, then for every $s\in F$
the inequality
\[
 s\wedge r\le s\ominus a
\]
holds in each chain factor. In fact, either the image of $r$ is
zero, or the image of $a$ is zero. The inequality follows in
both cases and hence in $A$. Since $s\wedge r\in F$, upward
closure gives $s\ominus a\in F$. Thus $D_F\subseteq N_F$.

We prove the quotient formula in (3). Let $I\subseteq N_F$ and suppose $q_I(a)\in D_{F_I}$. Choose $r\in F$ such that $q_I(a)\wedge q_I(r)=0$, and put
\[
 b=a\wedge r\in I,\qquad
 h=a\ominus r,\qquad
 d=r\ominus b=r\ominus a.
\]
Since $b\in I\subseteq N_F$ and $r\in F$, we have $d\in F$.
In the group representation, $h=(a-r)^+$ and $d=(r-a)^+$,
so $h\wedge d=0$ and $h\in D_F$. The identity
\[
 a=(a\wedge r)\oplus(a\ominus r)=b\oplus h
\]
then gives $a\in I\vee D_F$. Conversely, if $a\in D_F$, choose
$r\in F$ with $a\wedge r=0$. Applying $q_I$ shows that
$q_I(a)\wedge q_I(r)=0$, so $q_I(a)\in D_{F_I}$.
The inverse image $q_I^{-1}(D_{F_I})$ is an MV-ideal containing
both $D_F$ and $I$. Thus it contains $D_F\vee I$, and the two
inclusions give
\[
 q_I^{-1}(D_{F_I})=I\vee D_F.
\]
This proves the first formula. The second formula is Theorem~\ref{thm:intermediate-quotients}.

Taking $I=D_F$ gives $D_{p(F)}=\{0\}$, which says precisely that every image radius is a weak order unit. The remaining assertions of (1) follow from Theorem~\ref{thm:intermediate-quotients}. If $a\in D_F$, choose $r\in F$ with $a\wedge r=0$. Then $f(a)\wedge f(r)=0$, and the hypothesis on $f(r)$ forces $f(a)=0$. This proves (2). Finally, for $I\subseteq N_F$ the first quotient formula vanishes exactly when $D_F\subseteq I$, whereas the second vanishes exactly when $I=N_F$. The Hausdorff criterion completes the proof.
\end{proof}

The reflection is functorial. If
$f:(A,\tau_F)\longrightarrow(B,\tau_{F'})$ is a continuous
MV-homomorphism, then
$f(N_F)\subseteq N_{F'}$. Thus
\[
 A/N_F\longrightarrow B/N_{F'},\qquad
 q_A(a)\longmapsto q_B(f(a)),
\]
is a continuous MV-homomorphism. We next consider compact filter topologies.

\begin{corollary}\label{cor:compact}
Let $F$ be a proper lattice filter of an MV-algebra $A$. If
$(A,\tau_F)$ is compact, then $(A,\tau_F)$ is a topological
MV-algebra.
\end{corollary}

\begin{proof}
The Hausdorff reflection $q:A\longrightarrow A/N_F$ is a continuous
surjection, so its target is compact Hausdorff.
By \cite[Theorem~5.5]{LuanWeberYang2021}, the MV-operations on
$A/N_F$ are continuous. Since $q$ is a homomorphism,
\[
 q\circ\oplus_A=\oplus_{A/N_F}\circ(q\times q),
 \qquad q\circ(*_A)=(*_{A/N_F})\circ q
\]
are continuous. To transfer continuity back to $A$, let
$U\subseteq A$ be open. By Theorem~\ref{thm:reflection}, there is
an open set $V\subseteq A/N_F$ with $U=q^{-1}(V)$. Then
\[
 \oplus_A^{-1}(U)
 =(q\times q)^{-1}\bigl(\oplus_{A/N_F}^{-1}(V)\bigr)
\]
is open in $A\times A$, and
\[
 (*_A)^{-1}(U)=q^{-1}\bigl((*_{A/N_F})^{-1}(V)\bigr)
\]
is open in $A$. Thus MV-addition is jointly continuous and
negation is continuous.
\end{proof}

\begin{remark}\label{rem:compact-no-hausdorff}
Corollary~\ref{cor:compact} extends the continuity conclusion of
\cite[Theorem~5.5]{LuanWeberYang2021} to compact filter topologies
without assuming Hausdorffness. Indeed, the canonical quotient
$A/N_F$ is Hausdorff, and the inverse-image description of $\tau_F$
allows the continuity of the MV-operations to be lifted to $A$.
\end{remark}

For proper MV-filters, the basic neighborhoods are congruence classes
\cite{Asadzadeh2020}. The same holds for principal lattice filters, with
the polar of the generator as the congruence kernel.
\begin{proposition}\label{prop:principal}
Let $A$ be an MV-algebra, and let $0<f\in A$ and $F=\uparrow f$.
Then the following statements hold.
\begin{enumerate}
\item $N_F=\{a\in A:a\wedge f=0\}$ and, for every $x\in A$,
$B_F(x,f)=\{y\in A:d(x,y)\in N_F\}$.
\item $\tau_F$ is the partition topology of the congruence modulo
$N_F$, and its Hausdorff quotient is discrete.
\item $\tau_F$ is Hausdorff if and only if $a\wedge f=0$ implies
$a=0$. In this case $\tau_F$ is discrete.
\end{enumerate}
\end{proposition}
\begin{proof}
Since $f$ is the least radius, $a\in N_F$ if and only if
$f\ominus a\in F$, or equivalently $f\ominus a=f$.
In a chain factor with $f=0$, both this equality and
$a\wedge f=0$ hold. In a factor with $f>0$, the equality
$(f-a)^+=f$ forces $a=0$, and conversely $a=0$ implies the
equality. The subdirect representation therefore gives
$f\ominus a=f$ if and only if $a\wedge f=0$ in $A$.
Applying this equivalence to $d(x,y)$ yields
\[
 B_F(x,f)=\{y:d(x,y)\wedge f=0\}
         =\{y:d(x,y)\in N_F\}.
\]
This ball is the congruence class of $x$ modulo $N_F$ and is
contained in every other basic neighborhood of $x$. Thus the
topology is exactly the partition topology of these classes,
and the quotient is discrete. Finally, Hausdorffness is equivalent
to $N_F=\{0\}$, in which case every class is a singleton.
\end{proof}

For example, in $A=[0,1]^2$ with $f=(1,0)$, the first coordinate has the discrete topology and the second coordinate is indiscrete. This is the ideal-type topology considered in \cite{Asadzadeh2020,GanLuanDengYang2024}. Its corresponding ideal is the polar of $f$.

\section{Polar and relative infinitesimal obstructions}\label{sec:structure}
Section~3 produced two canonical ideals $D_F\subseteq N_F$. The ideal $N_F$ controls Hausdorff separation, whereas $D_F$ is the least kernel that turns all filter radii into weak order units. We now pass to the Mundici representation $A=\Gamma(G,u)$ to identify the group-theoretic counterparts of these ideals and to relate them to the ambient $C$-topology. Throughout this section, $G$ is an abelian lattice-ordered group with strong unit $u$; no divisibility assumption is imposed.

An \emph{admissible radius filter} is a nonempty lattice filter
$C\subseteq G^+$ such that $0\notin C$ and
\[
\text{for every }c\in C\text{ there exists }d\in C
\text{ with }2d\le c.
\]
By iteration, for every $m\ge1$ and $c\in C$ there exists
$d\in C$ such that $md\le c$.
If $G$ is 2-divisible, the condition is equivalent, by upward closure,
to closure of $C$ under halving.

For $x\in G$ and $c\in C$, put
\[
B_C^G(x,c)=\{y\in G:c-|x-y|\in C\}.
\]
The positive-filter topology of
\cite[Theorems~2.4 and~2.5]{KoppermanPajooheshRichmond2011}
has neighborhoods $\{y:n|x-y|\le c\}$, where $c\in C$ and $n\ge1$.
By the radius condition, the sets $\{y:|x-y|\le d\}$, $d\in C$,
form a neighborhood base. If $d\in C$ satisfies $2d\le c$, then
\[
\{g:|g|\le d\}\subseteq B_C^G(0,c)
\subseteq \{g:|g|\le c\}.
\]
If $y\in B_C^G(x,c)$, choose $d\in C$ such that
\[
2d\le c-|x-y|.
\]
Whenever $|z-y|\le d$, the triangle inequality yields
\[
c-|x-z|
\ge c-|x-y|-d
\ge d\in C.
\]
Hence these balls form an open base for the positive-filter topology,
which we denote by $\tau_C^G$. Under the hypotheses of
\cite{Gusic1998,Yang2009}, this is the corresponding $C$-topology.

Put $A=\Gamma(G,u)$ and $F=C\cap[0,u]$. The filter $F$ is
nonempty. Indeed, for
$c\in C$, choose $n\ge1$ with $c\le nu$ and then $d\in C$ with
$nd\le c$. Since multiplication by $n$ reflects the order \cite{AndersonFeil1988}, $d\le u$,
and hence $d\in F$.
Thus $F$ is a proper lattice filter of $A$. Fix $f\in F$. For every
$c\in C$, we have $c\wedge f\in F$ and $c\wedge f\le c$, so $F$ is
downward cofinal in $C$. In particular, $u\in F$.

For $x\in A$ and $r\in F$, the two ball definitions give
\[
B_C^G(x,r)\cap A\subseteq B_F(x,r)
=\{y\in A:(r-|x-y|)^+\in C\}.
\]
Consequently, $\tau_F\subseteq\tau_C^G|_A.$

To describe $N_F$ in the representing group, put
\[\mathcal N_C=\{g\in G:(c-|g|)^+\in C\text{ for every }c\in C\}.\]

\begin{lemma}\label{lem:group-kernel}
The set $\mathcal N_C$ is a proper $\ell$-ideal, and
$N_F=\mathcal N_C\cap[0,u].$   Consequently, by the standard ideal correspondence under the Mundici equivalence,
 $A/N_F\cong\Gamma(G/\mathcal N_C,u+\mathcal N_C).
$
\end{lemma}
\begin{proof}
We first verify that $\mathcal N_C$ is a solid subgroup.
The definition gives $0\in\mathcal N_C$ and
$g\in\mathcal N_C$ if and only if $-g\in\mathcal N_C$.
If $g\in\mathcal N_C$ and $|h|\le|g|$, then for every $c\in C$,
\[
 (c-|h|)^+\ge(c-|g|)^+\in C.
\]
Upward closure of $C$ gives $h\in\mathcal N_C$.
For closure under addition, let $g,h\in\mathcal N_C$ and $c\in C$.
Then $c_g=(c-|g|)^+\in C$, so
$(c_g-|h|)^+\in C$. The triangle inequality gives
\[
 (c-|g+h|)^+\ge(c-|g|-|h|)^+
 =\bigl((c-|g|)^+-|h|\bigr)^+
 =(c_g-|h|)^+\in C.
\]
Thus $\mathcal N_C$ is an $\ell$-ideal. It is disjoint from $C$ because
$c\in C\cap\mathcal N_C$ would imply $(c-|c|)^+=0\in C$.
Since $C\ne\varnothing$, this also proves that $\mathcal N_C$ is proper.

For $a\in[0,u]$, membership in $\mathcal N_C$ implies
$(r-a)^+\in C$ for every $r\in F$. Since
$0\le(r-a)^+\le r\le u$, this element belongs to $F$, hence $a\in N_F$.
Conversely, if $a\in N_F$ and $c\in C$, choose $r\in F$ with
$r\le c$. Then $(c-a)^+\ge(r-a)^+\in F\subseteq C$, so
$a\in\mathcal N_C$. The quotient identification now follows from
Proposition~\ref{prop:ideal-correspondence}.
\end{proof}

Using the principal polars from Section~\ref{sec:prelim}, define
$D_C=\bigcup_{c\in C}c^\perp.$
These polars form a directed family, since
$c^\perp,d^\perp\subseteq(c\wedge d)^\perp$. Their union $D_C$ is
therefore an $\ell$-ideal. We call it the \emph{polar obstruction
ideal}, although $D_C$ need not itself be a polar.

\begin{theorem}\label{thm:polar}
Let $p:G\longrightarrow H:=G/D_C$ and let $v=p(u)$  and  $E=p(C)$. Then the following statements hold.
\begin{enumerate}[label=(\arabic*)]
\item $D_C\subseteq\mathcal N_C$, and $E$ is an admissible radius filter in the unital $\ell$-group $(H,v)$.
\item Every element of $E$ is a weak order unit. Moreover, $D_C$ is the least $\ell$-ideal $I$ such that the image of every $c\in C$ is a nonzero weak order unit in $G/I$.
\item For every $g\in G$,
$g^+\in C$ if and only if $p(g)\in E.$
\end{enumerate}
\end{theorem}
\begin{proof}
(1) Let $g\in D_C$. Choose $c\in C$ with $|g|\wedge c=0$.
For $r\in C$, put $s=r\wedge c\in C$. The elements $s$ and
$r\wedge|g|$ are disjoint, so
\[
 s+(r\wedge|g|)=s\vee(r\wedge|g|)\le r.
\]
Consequently,
\[
 s\le r-(r\wedge|g|)=(r-|g|)^+.
\]
Upward closure gives $(r-|g|)^+\in C$. Hence
$D_C\subseteq\mathcal N_C$. By the proof of
Lemma~\ref{lem:group-kernel}, $C\cap\mathcal N_C=\varnothing$.
Hence $C\cap D_C=\varnothing$ and $0\notin p(C)=E$.

The set $E$ is nonempty and closed under finite meets. If
$p(c)\le h\in H^+$ with $c\in C$, choose $g\in G$ with $p(g)=h$.
Then $c\vee g\in C$ and $p(c\vee g)=h$, so $E$ is upward closed.
For $c\in C$, choose $d\in C$ with $2d\le c$. Then
$2p(d)\le p(c)$ and $p(d)\in E$, so $E$ satisfies the radius
condition. Also $v=p(u)$ is a strong unit. Thus $E$ is admissible.

\smallskip\noindent
(2) Let $c\in C$ and suppose that
$|p(g)|\wedge p(c)=0$. Then
\[
 p(|g|\wedge c)=|p(g)|\wedge p(c)=0,
\]
so $|g|\wedge c\in D_C$. Thus for some
$d\in C$,
\[
 |g|\wedge(c\wedge d)=0.
\]
As $c\wedge d\in C$, it follows that $g\in D_C$ and $p(g)=0$. Hence $p(c)$ is
a weak order unit. To prove minimality, let $I$ be an $\ell$-ideal such that
$\pi_I(c)=c+I$ is a nonzero weak order unit in $G/I$
for every $c\in C$, where $\pi_I:G\longrightarrow G/I$ is the quotient map.
For $g\in D_C$, choose $c\in C$ with $|g|\wedge c=0$. Then
\[
 |g+I|\wedge(c+I)=\pi_I(|g|\wedge c)=0,
\]
and therefore $\pi_I(g)=g+I=0$, so $g\in I$. Thus $D_C\subseteq I$.

\smallskip\noindent
(3) We first show that $p^{-1}(E)\cap G^+=C$. One inclusion
is immediate. For the other, let $b\in G^+$ and $p(b)=p(c)$ for
some $c\in C$. Then $b-c\in D_C\subseteq\mathcal N_C$, whence
\[
 (c-|b-c|)^+\in C,\qquad (c-|b-c|)^+\le b.
\]
Thus $b\in C$.

If $g^+\in C$, then $p(g^+)\in E$ is a weak order unit and
\[
 p(g^+)\wedge p(g^-)=p(g^+\wedge g^-)=0.
\]
Hence $p(g^-)=0$ and $p(g)=p(g^+)\in E$. Conversely, if
$p(g)\in E$, then
\[
 p(g^+)=p(g)^+=p(g)\in E,
 \qquad g^+\in p^{-1}(E)\cap G^+=C.
\]
This completes the proof.
\end{proof}

The polar quotient also describes the filter topology. The restriction
$p_A:A\longrightarrow\Gamma(H,v)$ is onto by
Proposition~\ref{prop:ideal-correspondence}, and $p(F)=E\cap[0,v]$.
Indeed, one inclusion follows from preservation of order. For the
other, if $p(c)\le v$ with $c\in C$, then $u\in C$ implies
$c\wedge u\in F$, and $p(c\wedge u)=p(c)\wedge v=p(c)$.
For $x,y\in A$ and $r\in F$, Theorem~\ref{thm:polar}(3) gives
\[
\begin{aligned}
 y\in B_F(x,r)
 &\Longleftrightarrow (r-|x-y|)^+\in C\\
 &\Longleftrightarrow p(r)-|p(x)-p(y)|\in E\\
 &\Longleftrightarrow p(y)\in B_E^H(p(x),p(r))\cap[0,v].
\end{aligned}
\]
Hence
\[
 B_F(x,r)=p_A^{-1}\bigl(B_E^H(p(x),p(r))\cap[0,v]\bigr).
\]
The filter $E\cap[0,v]$ is downward cofinal in $E$, by the argument
above applied to $(H,v,E)$. These balls therefore form a base of
the ambient subspace topology on $\Gamma(H,v)$. Thus $\tau_F$ is the inverse-image topology under $p_A$.
The MV-operations on $\Gamma(H,v)$ are continuous because they
are restrictions of the continuous group and lattice operations.
For an open set $V$ in this unit interval,
\[
 \oplus_A^{-1}\bigl(p_A^{-1}(V)\bigr)
 =(p_A\times p_A)^{-1}\bigl(\oplus_{\Gamma(H,v)}^{-1}(V)\bigr),
\]
and the corresponding identity holds for negation. It follows that
$A$ is a topological MV-algebra. For admissible radii this agrees
with the continuity criterion in
\cite[Theorem~3.5]{LuanWeberYang2021}. We do not assert this
conclusion for arbitrary lattice filters.

For an admissible radius filter $E$ in a unital $\ell$-group $H$, define
\[
 \operatorname{Inf}_E(H)=
 \{h\in H:n|h|\le e\text{ for all }e\in E\text{ and all integers }n\ge1\}.
\]
We call these elements infinitesimals relative to $E$. The
neighborhood description in
\cite[Theorems~2.4 and 2.5]{KoppermanPajooheshRichmond2011}
identifies this set with the intersection of all zero neighborhoods
in the positive-filter topology. It is consequently an $\ell$-ideal,
namely the class of $0$ under topological indistinguishability.
It may be a proper subideal of the usual infinitesimal ideal
determined by the unit.

\begin{theorem}\label{thm:decomposition}
Let $p:G\to H:=G/D_C$ and let $E=p(C)$ be as in
Theorem~\ref{thm:polar}. Then
$
\mathcal N_C
=
p^{-1}\bigl(\operatorname{Inf}_E(H)\bigr).
$
Consequently,
$
N_F
=
[0,u]\cap
p^{-1}\bigl(\operatorname{Inf}_E(H)\bigr).
$
Moreover, the restriction of $p$ yields the short exact sequence of
$\ell$-groups
\[
0\longrightarrow D_C
\xrightarrow{\ \iota\ }
\mathcal N_C
\xrightarrow{\ p|_{\mathcal N_C}\ }
\operatorname{Inf}_E(H)
\longrightarrow0,
\]
where $\iota$ is the inclusion map.
\end{theorem}
\begin{proof}
For $g\in G$, apply Theorem~\ref{thm:polar}(3) to $c-|g|$.
Since $p$ preserves absolute values and $E=p(C)$, we obtain
\[
 \begin{aligned}
 g\in\mathcal N_C
 &\Longleftrightarrow (c-|g|)^+\in C\quad(c\in C)\\
 &\Longleftrightarrow p(c)-|p(g)|\in E\quad(c\in C)\\
 &\Longleftrightarrow e-|p(g)|\in E\quad(e\in E).
 \end{aligned}
\]
For $h\in H$, the latter condition is equivalent to
$h\in\operatorname{Inf}_E(H)$. If $e-|h|\in E$ for every
$e\in E$, then for prescribed $e\in E$ and $n\ge1$ choose
$f\in E$ with $nf\le e$. Since $f-|h|\ge0$, one has
$n|h|\le nf\le e$. Conversely, if $h\in\operatorname{Inf}_E(H)$,
choose $f\in E$ with $2f\le e$. Then $|h|\le f$, and hence
$e-|h|\ge e-f\ge f\in E$. Upward closure proves the equivalence.

Applying it to $h=p(g)$ gives
$\mathcal N_C=p^{-1}(\operatorname{Inf}_E(H))$.
Intersecting with $[0,u]$ gives the asserted formula for $N_F$.
Since $\ker p=D_C\subseteq\mathcal N_C$, the restricted map
has kernel $D_C$. If $h\in\operatorname{Inf}_E(H)$, choose
$g\in G$ with $p(g)=h$. The inverse-image formula just proved
implies $g\in\mathcal N_C$. Thus the restricted map is onto
$\operatorname{Inf}_E(H)$, proving exactness.
\end{proof}

Consequently, the Hausdorff reflection factors through the polar quotient,
\[
 A\xrightarrow{\ p_A\ }\Gamma(H,v)
 \xrightarrow{\ q_E\ }
 \Gamma\bigl(H/\operatorname{Inf}_E(H),
              v+\operatorname{Inf}_E(H)\bigr)
 \cong A/N_F,
\]
where $q_E$ is the quotient homomorphism on the unit interval.
The first map preserves basic balls by Theorem~\ref{thm:polar},
and $q_E$ is the Hausdorff reflection for $E\cap[0,v]$ by
Theorem~\ref{thm:reflection}. Since $v\in E$,
$\operatorname{Inf}_E(H)\subseteq\operatorname{Inf}_v(H)$,
with equality whenever every element of $E$ is a strong unit.
Indeed, if $h\in\operatorname{Inf}_v(H)$ and $e\in E$ is a strong
unit, choose an integer $m\ge1$ with $v\le me$. For every $n\ge1$,
$mn|h|\le v\le me$, whence $n|h|\le e$.

We can now answer the subspace-topology question of
Luan, Zhao and Yang \cite[Question~2]{LuanZhaoYang2021}
under the radius refinement condition.
\Needspace{14\baselineskip}
\begin{theorem}\label{thm:weak-unit}
The following conditions are equivalent.
\begin{enumerate}[label=(\arabic*)]
\item $\tau_F=\tau_C^G|_{[0,u]}$.
\item Every element of $C$ is a weak order unit of $G$.
\item $D_C=\{0\}$.
\item For every $x\in[0,u]$ and $r\in F$,
$
 B_F(x,r)=B_C^G(x,r)\cap[0,u].
$
\end{enumerate}
Furthermore, $\tau_F$ is Hausdorff if and only if
\[D_C=\{0\}\quad\text{and}\quad\operatorname{Inf}_C(G)=\{0\}.\]
If $G$ is Archimedean, the four conditions above are also equivalent to Hausdorffness of $\tau_F$.
\end{theorem}
\begin{proof}
\smallskip
\noindent
\emph{$(2)\Longleftrightarrow(3)$.}
Recall that $D_C=\bigcup_{c\in C}c^\perp$, where
$c^\perp=\{g\in G:|g|\wedge c=0\}$. Since $c$ is a weak order unit
precisely when $c^\perp=\{0\}$, condition $(2)$ is equivalent to
$D_C=\{0\}$.

\smallskip
\noindent
\emph{$(3)\Longrightarrow(4)\Longrightarrow(1)$.}
Assume $D_C=\{0\}$. Then the quotient map in
Theorem~\ref{thm:polar} is the identity and $E=C$. Hence
Theorem~\ref{thm:polar}(3) gives
$g^+\in C$ if and only if $g\in C$. Taking
$g=r-|x-y|$ therefore yields
\[
B_F(x,r)=B_C^G(x,r)\cap[0,u],
\]
so $(4)$ holds.

We already know that $\tau_F\subseteq\tau_C^G|_{[0,u]}$. Conversely,
given $c\in C$, choose $r\in F$ with $r\le c$, using the downward
cofinality of $F$ in $C$. Then
\[
B_F(x,r)=B_C^G(x,r)\cap[0,u]
\subseteq B_C^G(x,c)\cap[0,u].
\]
Thus every basic neighborhood of the subspace topology contains a
$\tau_F$-neighborhood of the same point, and hence
$\tau_C^G|_{[0,u]}\subseteq\tau_F$. This proves $(1)$.

\smallskip
\noindent
\emph{$(1)\Longrightarrow(3)$.}
We argue by contraposition. Suppose $D_C\ne\{0\}$. Then there are
$c\in C$ and $0<g\in c^\perp$. Put $a=g\wedge u$. Since $u$ is a
strong unit, $a>0$. Otherwise $g\wedge u=0$, and hence
$g\wedge nu=0$ for every $n$, contradicting $0<g\le nu$ for some
$n$.

Now $0<a\le u$ and $a\perp c$, so $a\in D_C\subseteq\mathcal N_C$.
By Lemma~\ref{lem:group-kernel},
$N_F=\mathcal N_C\cap[0,u]$, and therefore $a\in N_F$. Hence $a$
belongs to every $\tau_F$-neighborhood of $0$.

However, $a\notin B_C^G(0,c)\cap[0,u]$. Indeed, otherwise
$c-a\in C\subseteq G^+$, so $a\le c$. Since $a\perp c$, this
would give $a=a\wedge c=0$, a contradiction. Thus the ambient
subspace topology has a neighborhood of $0$ not containing $a$,
whereas every $\tau_F$-neighborhood of $0$ contains $a$. Hence
$\tau_F\ne\tau_C^G|_{[0,u]}$, proving the contrapositive.

Thus conditions $(1)$ to $(4)$ are equivalent.

\smallskip
\noindent
\emph{Hausdorffness.}
We know that $\tau_F$ is Hausdorff if and only if $N_F=\{0\}$.
Since $N_F=\mathcal N_C\cap[0,u]$ and $u$ is a strong unit, every
nonzero $\ell$-ideal of $G$ meets $[0,u]$ nontrivially. Hence
$N_F=\{0\}$ if and only if $\mathcal N_C=\{0\}$.

Theorem~\ref{thm:decomposition} gives
$\mathcal N_C=\{0\}$ exactly when $D_C=\{0\}$ and
$\operatorname{Inf}_E(H)=\{0\}$. When $D_C=\{0\}$, we have
$H=G$ and $E=C$, which gives the stated Hausdorff criterion.

Finally, if $G$ is Archimedean, then $u\in C$ implies that every
$g\in\operatorname{Inf}_C(G)$ satisfies $n|g|\le u$ for all $n$.
Hence $g=0$, so $\operatorname{Inf}_C(G)=\{0\}$. Thus, in the
Archimedean case, Hausdorffness is equivalent to $D_C=\{0\}$, and
therefore to each of conditions $(1)$ to $(4)$.
\end{proof}

We next consider the quotients that preserve the filter. Since $F$
is downward cofinal in $C$, the ideal correspondence gives
$D_F=D_C\cap[0,u]$. Theorem~\ref{thm:mv-polar-quotient} therefore
implies the quotient formulas below. We give a group proof to make
the induced exact sequence explicit.

\begin{theorem}\label{thm:group-intermediate}
Let $I\subseteq\mathcal N_C$ be an $\ell$-ideal,
$\pi:G\longrightarrow\overline G=G/I$, and
$\overline C=\pi(C)$. Then $\overline C$ is an admissible radius
filter and
\[
 \pi^{-1}(\overline C)\cap G^+=C,\qquad
 \mathcal N_{\overline C}=\mathcal N_C/I,\qquad
 D_{\overline C}=(D_C+I)/I.
\]
The decomposition for $(\overline G,\overline C)$ is canonically
the exact sequence
\[
 0\longrightarrow D_C/(D_C\cap I)
 \longrightarrow\mathcal N_C/I
 \longrightarrow\mathcal N_C/(D_C+I)\longrightarrow0.
\]
Consequently, its unit-interval filter topology agrees with the
ambient subspace topology if and only if $D_C\subseteq I$.
It is Hausdorff if and only if $I=\mathcal N_C$.
\end{theorem}
\begin{proof}
If $b\in G^+$ and $\pi(b)=\pi(c)$ with $c\in C$, then
$b-c\in I\subseteq\mathcal N_C$ and
\[
 (c-|b-c|)^+\in C,\qquad (c-|b-c|)^+\le b.
\]
Hence $b\in C$, proving saturation. The image $\overline C$ is
nonempty and closed under finite meets. If
$\pi(c)\le\overline g\in\overline G^+$ with $c\in C$, choose
$g\in G$ with $\pi(g)=\overline g$. Then $c\vee g\in C$ and
$\pi(c\vee g)=\overline g$, proving upward closure.
Moreover, $I\cap C=\varnothing$ because
$I\subseteq\mathcal N_C$, so $0\notin\overline C$.
For each $c\in C$, choose $d\in C$ with $2d\le c$.
Then $2\pi(d)\le\pi(c)$, giving the radius condition in the quotient.

For $g\in G$, the definition of $\mathcal N_{\overline C}$ and
saturation applied to the nonnegative element $(c-|g|)^+$ give
\[
\begin{aligned}
 \pi(g)\in\mathcal N_{\overline C}
 &\Longleftrightarrow (\pi(c)-|\pi(g)|)^+\in\overline C
       \quad(c\in C)\\
 &\Longleftrightarrow \pi\bigl((c-|g|)^+\bigr)\in\overline C
       \quad(c\in C)\\
 &\Longleftrightarrow (c-|g|)^+\in C\quad(c\in C)\\
 &\Longleftrightarrow g\in\mathcal N_C.
\end{aligned}
\]
As $\pi$ is onto, this proves
$\mathcal N_{\overline C}=\mathcal N_C/I$.
To compute the polar obstruction ideal, suppose
$|\pi(g)|\wedge\pi(c)=0$ for some $c\in C$.
Put $b=|g|\wedge c\in I$. Since $b\in\mathcal N_C$,
\[
 d=c-b=(c-|g|)^+\in C,
 \qquad h=|g|-b=(|g|-c)^+\perp d.
\]
Thus $h\in D_C$ and $|g|=b+h\in I+D_C$, which implies
$g\in I+D_C$ by solidity of the sum of the two $\ell$-ideals.
Conversely, for $h\in D_C$ choose $c\in C$ with $|h|\wedge c=0$.
Applying $\pi$ gives $|\pi(h)|\wedge\pi(c)=0$, so
$\pi(h)\in D_{\overline C}$. Since $\pi(I)=\{0\}$, we obtain
$I+D_C\subseteq\pi^{-1}(D_{\overline C})$.
Thus $\pi^{-1}(D_{\overline C})=I+D_C$, proving the second ideal
formula.
Applying Theorem~\ref{thm:decomposition} to the quotient and
using the isomorphism theorems identifies its residual ideal with
$\mathcal N_C/(D_C+I)$ and its polar obstruction ideal with
$D_C/(D_C\cap I)$. More explicitly, the first map in the stated
sequence sends $d+(D_C\cap I)$ to $d+I$, and the second sends
$g+I$ to $g+(D_C+I)$. The first is injective, the second is onto,
and the kernel of the second is $(D_C+I)/I$, the image of the first.
Finally, $D_{\overline C}=\{0\}$ if and only if $D_C\subseteq I$,
so Theorem~\ref{thm:weak-unit} gives the subspace criterion.
The Hausdorff criterion is
$\mathcal N_{\overline C}=\{0\}$, which is equivalent to
$I=\mathcal N_C$ by the first ideal formula.
\end{proof}

\section{Concrete models for the obstruction ideals}\label{sec:models}

We compute the ideals in the preceding decomposition for finite
products of linearly ordered groups and for groups of continuous
functions. In the latter case, we also study when the associated
exact sequence splits.

\subsection{Finite products of linearly ordered groups}

Filter topologies on MV-chains are studied in
\cite[Section~4]{LuanWeberYang2021}. We classify admissible radius
filters on finite products of linearly ordered groups and determine
the associated ideals and topologies.

\begin{theorem}\label{thm:chain-products}
Let $G=\prod_{i=1}^nG_i$, where each $G_i$ is a nontrivial linearly
ordered abelian group with strong unit $u_i$, and set
$u=(u_1,\ldots,u_n)$. The admissible radius filters in $G^+$ are in
bijection with tuples of convex subgroups $K_i\le G_i$, not all
equal to $G_i$, such that every nonzero quotient $G_i/K_i$ has no
least positive element. For such a tuple, put $J=\{i:K_i\ne G_i\}$. The corresponding filter and ideals are
\begin{align*}
 C&=\{g\in G^+:g_i\notin K_i\text{ for every }i\in J\},\\
 D_C&=\prod_{i=1}^nD_i,\qquad
 D_i=\begin{cases}0&i\in J,\\G_i&i\notin J,\end{cases}\\
 \mathcal N_C&=\prod_{i=1}^nK_i.
\end{align*}
Let $A=\Gamma(G,u)$ and $F=C\cap A$. The Hausdorff quotient of $A$ is
\[A/N_F\cong\prod_{i\in J}\Gamma(G_i/K_i,u_i+K_i),\]
with the product of the order topologies. The topology $\tau_F$ is the inverse-image topology under the quotient map. Consequently,
\[
 \tau_F=\tau_C^G|_A\ \Longleftrightarrow\ J=\{1,\ldots,n\},
 \qquad
 \tau_F\text{ is Hausdorff}\ \Longleftrightarrow\ K_i=0\ (1\le i\le n).
\]
In the ambient topology, coordinates outside $J$ are discrete, while coordinates in $J$ carry the inverse images of the order topologies of $G_i/K_i$.
\end{theorem}
\begin{proof}
For each $i$, let $C_i$ be the $i$th coordinate projection of $C$.
This is a nonempty upward-closed subset of $G_i^+$. Indeed, if
$v_i\in C_i$ and $w_i\ge v_i$, replace the $i$th coordinate of
a vector $v\in C$ by $w_i$. The resulting vector lies above $v$
and belongs to $C$. The radius condition also passes to $C_i$.
For $c_i\in C_i$, choose $c\in C$ with this coordinate and
$d\in C$ with $2d\le c$. Then $d_i\in C_i$ and $2d_i\le c_i$.

First, we show that
\[
 C=\prod_{i=1}^n C_i.
\]
The inclusion $C\subseteq\prod_i C_i$ is clear. Conversely,
given $v_i\in C_i$ for every $i$, choose $c^{(i)}\in C$ with
$c^{(i)}_i=v_i$. The finite meet $c=\bigwedge_i c^{(i)}$ belongs
to $C$, and $c_i\le c^{(i)}_i=v_i$ for every $i$. Thus $c\le v$
for $v=(v_i)$, and upward closure gives $v\in C$.

If $0\in C_i$, put $K_i=G_i$, since then $C_i=G_i^+$.
Otherwise put
\[
 K_i=\{g\in G_i:|g|\notin C_i\}.
\]
The complement of $C_i$ in $G_i^+$ is downward closed, so $K_i$
is solid. To check addition, suppose that $a,b\ge0$ lie outside
$C_i$ but $a+b\in C_i$. Choose $d\in C_i$ with
\[
 2d\le a+b\le2\max\{a,b\}.
\]
Since multiplication by a positive integer reflects the order,
$d\le\max\{a,b\}$. Upward closure would give
$\max\{a,b\}\in C_i$, a contradiction. Thus $a+b\notin C_i$.
If $g,h\in K_i$, it follows that $|g|+|h|\notin C_i$, and
$|g+h|\le|g|+|h|$ gives $g+h\in K_i$. The definition also
gives $-g\in K_i$. Hence $K_i$ is a solid subgroup and therefore a convex subgroup. Thus
\[
 C_i=G_i^+\setminus K_i\quad(i\in J),
 \qquad J=\{i:K_i\ne G_i\},
\]
and the tuple is uniquely determined by $C$. Moreover, $J$ is
nonempty, since otherwise $0\in\prod_i C_i=C$.

For $i\in J$, let $\pi_i:G_i\longrightarrow G_i/K_i$ be the
quotient homomorphism. This quotient is linearly ordered by the
standard theory of convex subgroups \cite{AndersonFeil1988}.
It has no least positive element. If $\alpha$ were such an
element, choose $c_i\in C_i$ with $\pi_i(c_i)=\alpha$ and then
$d_i\in C_i$ with $2d_i\le c_i$. We would have
$0<\pi_i(d_i)$ and $2\pi_i(d_i)\le\alpha$, contradicting the
minimality of $\alpha$.

Conversely, suppose $(K_i)$ satisfies the stated conditions.
For $i\in J$, the set $G_i^+\setminus K_i$ is upward closed.
It is closed under finite meets because $G_i$ is linearly ordered.
Consequently, the displayed formula defines a lattice filter $C$.
It contains $u$, since a proper convex subgroup cannot contain
the strong unit $u_i$, and it does not contain $0$, since $J$
is nonempty. To verify the radius condition, fix $c\in C$.
For each $i\in J$, the absence of a least positive element in
$G_i/K_i$ allows us to choose
\[
 0<a_i<\pi_i(c_i),\qquad
 0<\delta_i<\min\{a_i,\pi_i(c_i)-a_i\}.
\]
Then $2\delta_i<\pi_i(c_i)$. Choose a lift $d_i$ of $\delta_i$.
In a quotient of a linearly ordered group by a convex subgroup,
a strictly positive quotient class has only positive
representatives \cite{AndersonFeil1988}. Applying this to $d_i$ and $c_i-2d_i$ gives
\[
 0<2d_i<c_i.
\]
Put $d_i=0$ for $i\notin J$. Then $d=(d_i)$ belongs to $C$
and satisfies $2d\le c$, as required.

For $i\in J$, every radius has
$c_i>0$, and the equality $|g_i|\wedge c_i=0$ in a chain forces
$g_i=0$. Hence every element of $D_C=\bigcup_{c\in C}c^\perp$
vanishes on $J$. Conversely, if $g_i=0$ for every $i\in J$,
choose $c\in C$ with $c_i=u_i$ for $i\in J$ and $c_i=0$
otherwise. Then $|g_i|\wedge c_i=0$ in every coordinate, so
$g\in c^\perp\subseteq D_C$. This gives the formula for $D_C$.

For $g\in G$ and $c\in C$, the condition
$(c-|g|)^+\in C$ is equivalent to
\[
 \pi_i\bigl((c_i-|g_i|)^+\bigr)>0\qquad(i\in J).
\]
Since each $\pi_i$ is an $\ell$-homomorphism, this is equivalent
to $(\pi_i(c_i)-|\pi_i(g_i)|)^+>0$, that is,
\[
 |\pi_i(g_i)|<\pi_i(c_i)\qquad(i\in J).
\]
If $g\in\mathcal N_C$ and $g_i\notin K_i$ for some $i\in J$,
choose a radius $c\in C$ with $c_i=|g_i|$, using
$C=\prod_i C_i$. The resulting inequality is
$|\pi_i(g_i)|<|\pi_i(g_i)|$, a contradiction. Therefore
$g_i\in K_i$ for all $i\in J$. Conversely, if $g_i\in K_i$ for every $i\in J$, then
$|\pi_i(g_i)|=0<\pi_i(c_i)$ for every $c\in C$.
Thus $g\in\mathcal N_C$. Since $K_i=G_i$ outside $J$, we obtain
$\mathcal N_C=\prod_iK_i$.

By the ideal correspondence in Proposition~\ref{prop:ideal-correspondence},
the quotient has the asserted form. Let $q$ denote its quotient map. For $x,y\in A$ and $r\in F$, the basic balls satisfy
\[
 \begin{aligned}
 y\in B_F(x,r)
 &\Longleftrightarrow (r-|x-y|)^+\in C\\
 &\Longleftrightarrow
 |\pi_i(x_i)-\pi_i(y_i)|<\pi_i(r_i)
 \quad\text{for every }i\in J.
 \end{aligned}
\]
Every quotient radius $0<\rho_i\le\pi_i(u_i)$ has a lift in
$[0,u_i]$. Indeed, truncate any representative to this interval,
using that $\pi_i$ preserves meets and joins. The truncated lift
still lies outside $K_i$. These lifts can be chosen independently,
and the other coordinates can be set to zero. Thus $q(F)$ consists
of all tuples of positive radii in the quotient unit intervals.
The displayed inequalities form the usual order-neighborhood bases
in each interval, including one-sided neighborhoods at its endpoints.
Since the product is finite, their products form a neighborhood
base for the product topology. The exact inverse-image identity
of Theorem~\ref{thm:reflection} therefore identifies the reflected
topology with this product and $\tau_F$ with its inverse image.

For the ambient topology, $c-|x-y|\in C$ is equivalent to the
same strict quotient inequalities for $i\in J$ together with
$|x_i-y_i|\le c_i$ for $i\notin J$. In the latter coordinates,
we can choose $c_i=0$, which forces $x_i=y_i$. The coordinates
of $c$ vary independently, so the ambient topology is the product
of the asserted inverse-image topologies and discrete topologies.
Finally, the formula for $D_C$ shows that $D_C=0$ exactly when
$J=\{1,\ldots,n\}$, and the formula for $\mathcal N_C$ shows
that $\mathcal N_C=0$ exactly when every $K_i=0$.
Theorem~\ref{thm:weak-unit} gives the subspace and Hausdorff criteria.
\end{proof}

The finiteness of the product is used in proving $C=\prod_i C_i$.
If $K_i=G_i$, the $i$th coordinate contributes to $D_C$. If
$0<K_i<G_i$, it contributes to $\mathcal N_C/D_C$.

Let $\mathbb R^m_{\mathrm{lex}}$ denote the additive group with lexicographic order. For $0\le k\le m$, put
\[
 K^{(k)}=\{x\in\mathbb R^m:x_1=\cdots=x_k=0\},
 \qquad K^{(0)}=\mathbb R^m,\quad K^{(m)}=0.
\]
These are precisely its convex subgroups \cite{AndersonFeil1988}.

\begin{corollary}\label{cor:lex}
Take $G=\prod_{i=1}^n\mathbb R^{m_i}_{\mathrm{lex}}$ and $u_i=(1,0,\ldots,0)$. Every admissible $C$ is uniquely specified by integers $0\le k_i\le m_i$, not all zero, through Theorem~\ref{thm:chain-products} with $K_i=K_i^{(k_i)}$. Its separation ideal is
\[
 N_F=[0,u]\cap\prod_iK_i^{(k_i)},
\]
and its Hausdorff quotient is
\[
 \prod_{k_i>0}\Gamma(\mathbb R^{k_i}_{\mathrm{lex}},(1,0,\ldots,0))
\]
with product order topology. The subspace equality holds exactly when every $k_i>0$. Hausdorffness holds exactly when every $k_i=m_i$. Compactness holds exactly when every $k_i\in\{0,1\}$.
\end{corollary}
\begin{proof}
The quotient description and the subspace and Hausdorff criteria
follow from Theorem~\ref{thm:chain-products}, since
$\mathbb R^{m_i}_{\mathrm{lex}}/K_i^{(k_i)}
\cong\mathbb R^{k_i}_{\mathrm{lex}}$ for $k_i>0$, and every
such quotient has no least positive element. In particular,
$K_i^{(k_i)}$ is proper exactly when $k_i>0$, and it is zero
exactly when $k_i=m_i$.
For compactness, suppose first that $k_i=1$ whenever $k_i>0$.
The Hausdorff quotient is then a finite product of copies of $[0,1]$
and is compact. Suppose now that $k_i\ge2$ for some $i$. Then
\[
\Gamma(\mathbb R_{\mathrm{lex}}^{k_i},(1,0,\ldots,0))
\]
is not compact in its order topology. For $t\in[0,1]$, let
$V_t$ consist of the points of this interval with first coordinate
$t$. Each $V_t$ is nonempty. If a point of $V_t$ is neither
endpoint of the unit interval, choose smaller and larger points
in the same fiber by changing a later coordinate. The open order
interval between them remains in $V_t$. At the two endpoints, use
one-sided order neighborhoods with the other endpoint in the same
fiber. Thus the sets $V_t$ form an infinite pairwise disjoint open
cover with no finite subcover. The product quotient is also
noncompact, since it maps continuously onto this factor.

Finally, $\tau_F$ is the inverse-image topology under the canonical
Hausdorff quotient map $q$. Compactness of $A$ implies compactness
of the quotient by continuity and surjectivity. Conversely, every
open cover of $A$ has the form $\{q^{-1}(U_\lambda)\}_\lambda$
for open sets $U_\lambda$ in the quotient. As $q$ is surjective, the sets $U_\lambda$ cover the quotient.
A finite subcover pulls back to a finite subcover of $A$. Thus
\[
(A,\tau_F)\text{ is compact}
\quad\Longleftrightarrow\quad
k_i\in\{0,1\}\text{ for every }i.
\]
\end{proof}
When every $m_i=1$, the admissible sets in $\mathbb R^n$ are
\[
 C_J=\{c\in\mathbb R_+^n:c_i>0\ (i\in J)\},\qquad\varnothing\ne J\subseteq\{1,\ldots,n\}.
\]
Here $N_F=\{a\in[0,1]^n:a_i=0\ (i\in J)\}$. The filter topology is the product of the usual topology in coordinates in $J$ and the indiscrete topology in the others. The ambient topology has discrete omitted coordinates, illustrating the discrepancy noted in \cite[Remark~3.1]{LuanZhaoYang2021}.

The following examples show that the Archimedean hypothesis in
Theorem~\ref{thm:weak-unit} cannot be omitted. Both are covered by
Theorem~\ref{thm:chain-products}. The first is adapted from
\cite[Example~6.1.7]{Yang2006}.

\begin{example}\label{ex:lex-nonhausdorff}
Let $G=\mathbb R^2_{\mathrm{lex}}$, $u=(1,0)$, and
\[
 C=\{(s,t):s>0\}.
\]
Here $K=\{0\}\times\mathbb R$ in Theorem~\ref{thm:chain-products}.
Every element of $C$ is a weak order unit, since $G$ is a chain. Nevertheless,
\[
 D_C=0,\qquad \mathcal N_C=\{0\}\times\mathbb R,\qquad
 N_F=\{(0,t):t\ge0\}.
\]
The two topologies on $A=\Gamma(G,u)$ coincide, but neither is Hausdorff. Both are the inverse images of the usual topology on $[0,1]$ under $(s,t)\mapsto s$. In particular, this is a compact, non-Hausdorff filter topology. For $a=(0,1)>0$, one has $na\le u$ for every positive integer $n$. This shows the failure of
\[
 \text{all admissible radii weak order units}
 \ \Longrightarrow\ \text{Hausdorffness}.
\]

\end{example}

\begin{example}\label{ex:both}
Let $G=\mathbb R^2_{\mathrm{lex}}\times\mathbb R$, $u=((1,0),1)$, and
$ C=\{((s,t),r):s>0,\ r\ge0\}.$
The tuple in Theorem~\ref{thm:chain-products} is
$K_1=\{0\}\times\mathbb R$ and $K_2=\mathbb R$, with
$J=\{1\}$. Therefore
\[
 D_C=\{(0,0)\}\times\mathbb R,
 \qquad
 \mathcal N_C=(\{0\}\times\mathbb R)\times\mathbb R.
\]
The last coordinate is the disjointness obstruction. After it is removed, the second lexicographic coordinate is the relative infinitesimal obstruction. Thus $0<D_C<\mathcal N_C$, the subspace equality fails, and $\tau_F$ is not Hausdorff. Its Hausdorff quotient is again the standard algebra $[0,1]$.
\end{example}

\begin{example}\label{ex:nondivisible}
The radius condition can hold in a nondivisible group. Take
\[
 G=\mathbb Z[1/3]\mathbin{\overrightarrow{\times}}\mathbb Z,
 \qquad u=(1,0),
 \qquad C=\{(a,b)\in G:a>0\},
\]
where the order is lexicographic. For $c=(a,b)\in C$, the element
$d=(a/3,0)$ belongs to $C$ and satisfies $2d<c$. The group $G$ is
not 2-divisible, since $(0,1)$ has no half. By
Theorem~\ref{thm:chain-products},
\[
 D_C=0,\qquad \mathcal N_C=\{0\}\times\mathbb Z.
\]
Thus a nonzero relative infinitesimal obstruction occurs outside the
2-divisible setting. The Hausdorff quotient of the unit interval is
$\Gamma(\mathbb Z[1/3],1)$ with its order topology. Here the
quotient group $G/(\{0\}\times\mathbb Z)\cong\mathbb Z[1/3]$
has no least positive element, as required by the classification.
\end{example}

\begin{remark}\label{rem:radius-condition-needed}
The condition on small radii cannot simply be dropped from the
relative infinitesimal formula. In
$G=\mathbb R\mathbin{\overrightarrow{\times}}\mathbb R$ with
$u=(1,0)$, let $C=\{c\in G:c\ge u\}$. Then $C$ is a proper
lattice filter, but its least element $u$ prevents it from satisfying
the radius condition. Since
$(u-|g|)^+<u$ whenever $g\ne0$, one has
\[
 D_C=\mathcal N_C=0,
 \qquad \operatorname{Inf}_C(G)=\{0\}\times\mathbb R.
\]
In particular, the separation ideal need not equal the relative
infinitesimal ideal in the absence of sufficiently small radii.
Indeed, for $g=(a,b)$ the inequalities $n|g|\le u$ for every
$n\ge1$ hold exactly when $a=0$. Since $u$ is the least radius,
these are precisely the inequalities defining
$\operatorname{Inf}_C(G)$.
\end{remark}

\subsection{Function groups and splitting}

We now study sections of the exact sequence for function groups.
A positive additive section of a quotient homomorphism
$q:N\longrightarrow R$ is an additive map $s:R\longrightarrow N$ with $q\circ s=\mathrm{id}_R$
and $s(R^+)\subseteq N^+$. Such a section need not preserve the
lattice operations.

\begin{theorem}\label{thm:function-positive-splitting}
Let $X$ be compact Hausdorff, let $G$ be a unital lattice-ordered
subring of $C(X,\mathbb R)$ with unit $\mathbf1$, and let $C$ satisfy
the hypotheses of Theorem~\ref{thm:decomposition}. Put
$D=D_C$, $N=\mathcal N_C$ and $R=N/D$.
If $0<a\in R$, there is $0<b\in R$ such that
\[
 nb\le a\qquad\text{for every integer }n\ge1.
\]
Consequently, the following conditions are equivalent.
\begin{enumerate}
\item the quotient map $N\longrightarrow R$ admits a positive
additive section.
\item the sequence $0\longrightarrow D\longrightarrow N
\longrightarrow R\longrightarrow0$ splits as a sequence of
$\ell$-groups.
\item $N=D$.
\item $R$ is Archimedean.
\end{enumerate}
More generally, $R$ admits an injective positive homomorphism into an
Archimedean $\ell$-group if and only if $R=0$.
When $G=C(X,\mathbb R)$, the sequence always has a real-linear section
and hence splits as a sequence of abelian groups.
\end{theorem}

\begin{proof}
Write $p:G\longrightarrow G/D$ and $v=p(\mathbf1)$.
Radius refinement and the strong-unit property give
$\mathbf1\in C$. Theorem~\ref{thm:decomposition} therefore yields
\[
 n|p(f)|\le v\qquad(f\in N,\ n\ge1).
\]
Choose $f_0\in N$ with $p(f_0)=a>0$. Since $p$ preserves
positive parts, $f=f_0^+\in N^+$ satisfies $p(f)=a$. Put
$b=p(f^2)$.
The square belongs to $G$ because $G$ is a subring. Choose a
positive integer $m$ with $f\le m\mathbf1$. Then $f^2\le mf$,
so $f^2\in N$.
If $f^2\in D$, then $f^2\perp c$ for some $c\in C$.
Since $f$ and $f^2$ have the same zero set, $f\perp c$ and hence
$f\in D$, contrary to $a>0$.
Thus $b>0$.

For each positive integer $n$, the inequality $na\le v$ gives
$(nf-\mathbf1)^+\in D$. The pointwise identity
\[
 (nf^2-f)^+=f(nf-\mathbf1)^+
 \le m(nf-\mathbf1)^+
\]
and solidity of $D$ imply $(nf^2-f)^+\in D$, proving $nb\le a$.
In particular, $R$ is non-Archimedean whenever it is nonzero.
If $\varphi:R\to L$ is a positive homomorphism and $L$ is
Archimedean, then, for every $n\ge1$,
\[
0\le n\varphi(b)\le\varphi(a).
\]
Hence $\varphi(b)=0$. Such a homomorphism cannot be injective when
$R\ne0$. Since $G$ is a lattice-ordered subgroup of the Archimedean
$\ell$-group $C(X,\mathbb R)$, it is Archimedean, and so is its
$\ell$-ideal $N$. In particular, a positive additive section of
$N\to R$ would be an injective positive homomorphism from $R$ into
the Archimedean group $N$. Thus (1) implies (3). If $N=D$, the zero
map is an $\ell$-group section, so (3) implies (2), and every
$\ell$-group section is positive, giving (2) implies (1).
The inequality $nb\le a$ also shows that (4) implies (3), while
the zero group is Archimedean. This proves the equivalences and the assertion about injective
positive homomorphisms.
Finally, for $G=C(X,\mathbb R)$, solidity makes both $D$ and $N$
real vector subspaces. For $f$ in either ideal and $\lambda\in\mathbb R$,
one has $|\lambda f|\le m|f|$ for some integer $m$.
A linear complement of $D$ in $N$ supplies the asserted section.
\end{proof}

\begin{remark}
The argument also shows that every nonzero $\ell$-ideal of $R$ is
non-Archimedean. If it contains $a>0$, it also contains the element
$0<b\le a$ constructed above.
\end{remark}

The multiplicative closure in
Theorem~\ref{thm:function-positive-splitting} cannot be omitted.

\begin{example}\label{ex:piecewise-affine-splitting}
Let $G$ be the real vector lattice of continuous piecewise-affine
functions on $[0,1]$ with finitely many pieces, let $u=\mathbf1$,
and put $C=\{c\in G^+:c(0)>0\}$.
The group $G$ is Archimedean and divisible, and $C$ is a proper
lattice filter closed under halving. We have
\[
 D_C=\{g\in G:g=0\text{ near }0\},
 \qquad
 \mathcal N_C=\{g\in G:g(0)=0\}.
\]
Indeed, if $g\perp c$ for $c\in C$, then $c$ is strictly positive
near $0$, so $g$ vanishes there. Conversely, if $g=0$ on
$[0,\delta]$ for some $\delta>0$, the function
$c(t)=(\delta-t)^+$ belongs to $C$ and is disjoint from $g$.
For the second identity, $(c-|g|)^+(0)>0$ for every $c\in C$
exactly when $g(0)=0$. If $g(0)\ne0$, the constant radius
$|g(0)|/2$ violates this condition.

Every $g\in G$ is affine on a sufficiently small interval beginning
at $0$. Hence
\[
 g+D_C\longmapsto\bigl(g(0),g'_+(0)\bigr)
\]
is an isomorphism $G/D_C\cong\mathbb R^2_{\mathrm{lex}}$. The map
is additive and surjective, since the function $a+bt$ realizes
the pair $(a,b)$. Its kernel consists exactly of functions
vanishing near $0$. Moreover, $g+D_C\ge0$ if and only if
$g^-\in D_C$, equivalently $g\ge0$ near $0$. For an affine
germ this holds exactly when $g(0)>0$, or $g(0)=0$ and
$g'_+(0)\ge0$. Thus the induced bijection and its inverse
preserve order, and hence preserve the lattice operations.
The residual ideal is
$\mathcal N_C/D_C\cong\{0\}\times\mathbb R$.
The quotient map on $\mathcal N_C$ is $g\mapsto g'_+(0)$, and
\[
 s:\mathbb R\longrightarrow\mathcal N_C,
 \qquad s(r)(t)=rt,
\]
is an $\ell$-group section, since taking a right derivative at
$0$ gives $r$ and multiplication by $t\ge0$ preserves the order
of the scalar coefficients. The group $G$ is not closed under
multiplication, as the function $t^2$ is not piecewise affine.
Thus a nonzero residual infinitesimal
ideal can admit a section even when the original group is
Archimedean.
\end{example}

Let $X$ be a compact Hausdorff space. The group $C(X,\mathbb R)$, with pointwise order and unit $\mathbf1$, is Archimedean. Its positive weak order units are precisely the functions whose cozero sets are dense in $X$. See \cite[Examples~3.1(b)(ii)]{KoppermanPajooheshRichmond2011}.
Hence Theorem~\ref{thm:weak-unit} gives, for any admissible $C\subseteq C(X,\mathbb R)^+$,
\[
 \tau_F=\tau_C^G|_{[0,\mathbf1]}
 \quad\Longleftrightarrow\quad
 \operatorname{coz}(c)\text{ is dense for every }c\in C
 \quad\Longleftrightarrow\quad \tau_F\text{ is Hausdorff}.
\]
For a closed set $K\subseteq X$, the following filter gives a
description of both ideals.

\begin{theorem}\label{thm:functions}
Let $X$ be a compact Hausdorff space and let
$\varnothing\neq K\subseteq X$ be closed. Put
$G=C(X,\mathbb R)$, $u=\mathbf 1$, and
\[
C_K=\{c\in C(X,\mathbb R)^+:c(k)>0\text{ for every }k\in K\},
\qquad
F_K=C_K\cap C(X,[0,1]).
\]
Then $C_K$ is admissible and
\[
\begin{aligned}
D_{C_K}
 &=\{g\in C(X,\mathbb R):g=0\text{ on some open neighborhood of }K\},\\
\mathcal N_{C_K}
 &=\{g\in C(X,\mathbb R):g|_K=0\}.
\end{aligned}
\]

The restriction homomorphism
\[
q:C(X,[0,1])\longrightarrow C(K,[0,1]),
\qquad q(f)=f|_K,
\]
induces an isomorphism
\[
C(X,[0,1])/N_{F_K}\cong C(K,[0,1]).
\]
Moreover, if $H=C(X,\mathbb R)/D_{C_K}$,
$p:C(X,\mathbb R)\longrightarrow H$ is the quotient map,
$v=p(\mathbf 1)$, and $E=p(C_K)$, then $H$ is the ordered
group of germs along $K$, every element of $E$ is a strong unit, and
\[
\mathcal N_{C_K}/D_{C_K}
\cong
\operatorname{Inf}_E(H)
=
\operatorname{Inf}_v(H).
\]
\end{theorem}

\begin{proof}
The admissibility of $C_K$ follows from the pointwise order, finite
meets, and halving. Moreover, $\mathbf 1\in C_K$, while
$0\notin C_K$ because $K\ne\varnothing$.

We first compute $D_{C_K}$. Suppose that $g\perp c$ for some
$c\in C_K$. Since $c>0$ on $K$, its cozero set
\[
\operatorname{coz}(c)=\{x\in X:c(x)\ne0\}
\]
is an open neighborhood of $K$. Since $|g|\wedge c=0$ pointwise, $g$ vanishes on
$\operatorname{coz}(c)$.

Conversely, suppose that $g=0$ on an open neighborhood $U$ of $K$.
The space $X$ is normal because it is compact Hausdorff.
By Urysohn's lemma \cite{Engelking1989}, there exists $c\in C(X,[0,1])$ such that
\[
c|_K=1,
\qquad
c|_{X\setminus U}=0.
\]
Then $c\in C_K$, and $\operatorname{coz}(c)\subseteq U$. Since $g$
vanishes on $U$, we have $|g|\wedge c=0$, so $g\perp c$.

If $g|_K=0$, then for every
$c\in C_K$,
\[
(c-|g|)^+|_K=c|_K>0.
\]
Hence $(c-|g|)^+\in C_K$ for every $c\in C_K$, and therefore
$g\in\mathcal N_{C_K}$.

Conversely, suppose that $g(k_0)\ne0$ for some $k_0\in K$. Choose
$0<\varepsilon<|g(k_0)|$ and put
$c=\varepsilon\mathbf 1$. Then $c\in C_K$, but
\[
(c-|g|)^+(k_0)=0.
\]
Thus $(c-|g|)^+\notin C_K$, and hence
$g\notin\mathcal N_{C_K}$.

The restriction map $q$ preserves the MV-operations because they
are defined pointwise. It is surjective by the Tietze extension
theorem \cite{Engelking1989}, applied with range $[0,1]$.
Lemma~\ref{lem:group-kernel} and the formula for
$\mathcal N_{C_K}$ give
\[
 \ker q=\{f\in C(X,[0,1]):f|_K=0\}
       =\mathcal N_{C_K}\cap[0,\mathbf1]=N_{F_K}.
\]
The first isomorphism theorem gives the asserted MV-algebra quotient.

Two functions have the same image in $H$ if and only if they agree
on a neighborhood of $K$. Thus $H$ is the group of germs along $K$. This identification preserves
the order, since
\[
 p(f)\le p(g)
 \quad\Longleftrightarrow\quad (f-g)^+\in D_{C_K}
 \quad\Longleftrightarrow\quad
 f\le g\text{ on some neighborhood of }K.
\]
Conversely, if a continuous function is defined on an open
neighborhood $U$ of $K$, choose a closed neighborhood $L$ of $K$
with $L\subseteq U$. Its restriction to $L$ extends to $X$ by the
Tietze extension theorem. Hence every such germ has a representative
in $C(X,\mathbb R)$.

Let $c\in C_K$. Since $K$ is compact and $c>0$ on $K$, there exists
$\varepsilon>0$ such that
\[
c(k)\ge 2\varepsilon
\qquad(k\in K).
\]
By continuity, after shrinking to a neighborhood $U$ of $K$ if
necessary, we have
$
c\ge\varepsilon\mathbf 1
$ on $U$.
Hence
$
(\varepsilon\mathbf 1-c)^+
$
vanishes on $U$ and therefore belongs to $D_{C_K}$. It follows that
$
p(c)\ge\varepsilon v
$
in $H$.

Since $\mathbf 1$ is a strong unit of $C(X,\mathbb R)$, its image
$v=p(\mathbf 1)$ is a strong unit of $H$. Therefore every positive
element dominating a positive scalar multiple of $v$ is also a strong
unit. Hence every element of
$E=p(C_K)$
is a strong unit.

By Theorem~\ref{thm:decomposition},
\[
\mathcal N_{C_K}/D_{C_K}
\cong
\operatorname{Inf}_E(H).
\]
Since every element of $E$ is a strong unit, the strong-unit case gives
$
\operatorname{Inf}_E(H)
=
\operatorname{Inf}_v(H).
$
\end{proof}

Under the restriction map, the image of $F_K$ consists of all
strictly positive elements of $C(K,[0,1])$. Indeed, each such
function has an extension to $C(X,[0,1])$, and every extension
belongs to $F_K$. This filter induces the uniform norm topology
on the compact space $K$, as in
\cite[Examples~3.1(b)(i) and (iii)]{KoppermanPajooheshRichmond2011}.
More explicitly, its ball with radius $r$ is
\[
 \{h\in C(K,[0,1]):|h(k)-f(k)|<r(k)\text{ for every }k\in K\}.
\]
The positive minimum of $r$ provides a uniform norm ball inside
this set, while positive constant radii give uniform norm balls.
Theorem~\ref{thm:reflection} therefore identifies $\tau_{F_K}$
with the inverse-image topology under $q$. In particular,
\[
\tau_{F_K}\text{ is Hausdorff}
\quad\Longleftrightarrow\quad
\tau_{F_K}=\tau_{C_K}^{G}|_{C(X,[0,1])}
\quad\Longleftrightarrow\quad
K=X.
\]
If $K=X$, every element of $C_K$ is strictly positive and hence
a weak order unit. If $K\ne X$, choose a nonempty open set $V$
whose closure is disjoint from $K$. Urysohn's lemma gives
$c\in C(X,[0,1])$ with $c|_K=1$ and $c|_{\overline V}=0$.
Then $c\in C_K$, but its cozero set is not dense. Thus
Theorem~\ref{thm:weak-unit} gives the displayed equivalences.

\begin{example}\label{ex:germs}
Take $X=[0,1]$, $K=\{0\}$, and $g(t)=t$. Then
$g\in\mathcal N_{C_K}$ but $g\notin D_{C_K}$, since $g$ vanishes
at $0$ but on no neighborhood of $0$. Thus $p(g)$ is a nonzero
positive infinitesimal in $H$. For every integer $n\ge1$,
$(ng-\mathbf1)^+$ vanishes on $[0,1/n]$, so $np(g)\le v$.
This shows that the relative infinitesimal term cannot be computed
inside the original Archimedean group before taking the polar
quotient. The Hausdorff reflection here is evaluation at $0$ onto
$[0,1]$.

Moreover, $D_{C_K}$ is not itself a polar. For every $t_0>0$,
the function $h_{t_0}(t)=(t-t_0/2)^+$ belongs to $D_{C_K}$ and
is nonzero at $t_0$. A continuous function disjoint from all
of $D_{C_K}$ therefore vanishes on $(0,1]$ and hence on $X$.
Thus $D_{C_K}^{\perp}=\{0\}$ and
$D_{C_K}^{\perp\perp}=G\ne D_{C_K}$, where
$S^\perp=\{g:g\perp s\text{ for all }s\in S\}$.
\end{example}

The next family has a fixed polar obstruction ideal but distinct separation ideals.

\begin{example}\label{ex:vanishing-rates}
Let $G=C([0,1],\mathbb R)$, with unit $\mathbf1$. For each
$\alpha>0$, put $w_\alpha(t)=t^\alpha$ and
\[
 C_\alpha=
 \{c\in G^+:\text{there exist }\varepsilon,\delta>0
 \text{ such that }c(t)\ge\varepsilon t^\alpha
 \text{ for }0<t<\min\{\delta,1\}\}.
\]
Each $C_\alpha$ contains $w_\alpha$ and excludes $0$. Upward
closure is immediate. If two radii satisfy the defining estimate
with parameters $(\varepsilon_1,\delta_1)$ and
$(\varepsilon_2,\delta_2)$, their meet satisfies it with
$\varepsilon=\min\{\varepsilon_1,\varepsilon_2\}$ and
$\delta=\min\{\delta_1,\delta_2\}$. Halving a radius replaces
$\varepsilon$ by $\varepsilon/2$. Thus these are admissible
radius filters. Their polar obstruction ideal is independent of $\alpha$ and is given by
\[
 D_{C_\alpha}=D
 =\{g\in G:g=0\text{ on a neighborhood of }0\}.
\]
Their separation ideals are
\[
 \mathcal N_{C_\alpha}=N_\alpha
 =\{g\in G:g(t)=o(t^\alpha)\text{ as }t\downarrow0\}.
\]
Indeed, if $g\perp c$ with $c\in C_\alpha$, then $g$ vanishes
on a punctured neighborhood of $0$ and hence, by continuity, on a
neighborhood of $0$. Conversely, if $g$ vanishes near $0$, choose a
continuous cutoff $\eta$ which is $1$ on a smaller neighborhood of
$0$ and whose support is contained in a neighborhood where $g=0$.
Then $\eta w_\alpha\in C_\alpha$ and
$g\perp\eta w_\alpha$.

If $g=o(t^\alpha)$ and $c(t)\ge\varepsilon t^\alpha$ near $0$,
then $(c-|g|)^+\ge(\varepsilon/2)t^\alpha$ on a sufficiently small
punctured neighborhood. Thus $g\in\mathcal N_{C_\alpha}$.
Conversely, apply the defining condition of $\mathcal N_{C_\alpha}$
to $c=\varepsilon w_\alpha$ for each $\varepsilon>0$.
Membership of $(\varepsilon w_\alpha-|g|)^+$ in $C_\alpha$ implies
$|g(t)|<\varepsilon t^\alpha$ for all sufficiently small $t>0$.
Hence $g=o(t^\alpha)$.

Let $p:G\longrightarrow H=G/D$, $v=p(\mathbf1)$,
$e_\alpha=p(w_\alpha)$ and $E_\alpha=p(C_\alpha)$.
Every element of $E_\alpha$ dominates a positive constant multiple
of $e_\alpha$, because the defining inequality holds near $0$.
All positive constant multiples of $e_\alpha$ belong to $E_\alpha$.
If $n|h|\le e_\alpha$ for every $n\ge1$, then for any
$\lambda>0$ and $n\ge1$ choose an integer $m\ge1/\lambda$.
The inequality $mn|h|\le e_\alpha$ gives
$n|h|\le\lambda e_\alpha$. Thus testing against $e_\alpha$
is equivalent to testing against all of $E_\alpha$.
Consequently,
\[
 \operatorname{Inf}_{E_\alpha}(H)
 =\{h\in H:n|h|\le e_\alpha\text{ for every }n\ge1\}
 =N_\alpha/D.
\]
Writing $N_0=\{g\in G:g(0)=0\}$, Theorem~\ref{thm:functions}
gives $\operatorname{Inf}_v(H)=N_0/D$. For $0<\alpha<\beta$,
\[
 D\subsetneq N_\beta\subsetneq N_\alpha\subsetneq N_0.
\]
The first, second and third strict inclusions are witnessed,
respectively, by $t^{\beta+1}$, $t^{(\alpha+\beta)/2}$ and
$t^\alpha$. Thus a single polar quotient supports a continuum of
distinct relative infinitesimal ideals, all strictly smaller than
its ordinary strong-unit infinitesimal ideal.
The associated exact sequences have linear sections but no positive
additive sections, by Theorem~\ref{thm:function-positive-splitting}.
These filters are not of the form $C_K$. Indeed, the functions $t^{2\alpha}$
and $t^\alpha$ have the same zero set, while only the first belongs
to $N_\alpha$.
\end{example}

For the continuous-function construction, splitting of the exact
sequence in Theorem~\ref{thm:decomposition} is characterized by
the following property of $K$ \cite[p.~294]{HartLuoVanMill2014}.

\begin{definition}\label{def:p-set}
A closed subset $K$ of a topological space $X$ is a \emph{$P$-set} if
the intersection of every countable family of neighborhoods of $K$ is
again a neighborhood of $K$.
\end{definition}

\Needspace{13\baselineskip}
\begin{proposition}\label{prop:germ-nonsplit}
Let $X$ be a compact Hausdorff space and $K\subseteq X$ a nonempty
closed set. Put
\[
 N=\{f\in C(X,\mathbb R):f|_K=0\},\qquad
 D=\{f\in C(X,\mathbb R):f=0\text{ on a neighborhood of }K\}.
\]
Then $N=\mathcal N_{C_K}$ and $D=D_{C_K}$. The following conditions
are equivalent.
\begin{enumerate}
\item The short exact sequence of $\ell$-groups
$0\longrightarrow D\longrightarrow N\longrightarrow N/D\longrightarrow0$
splits.
\item $N=D$.
\item $K$ is a $P$-set in $X$.
\item $C(X,\mathbb R)/D$ is Archimedean.
\end{enumerate}
As a sequence of abelian groups, the sequence always splits.
\end{proposition}
\begin{proof}
The identities $N=\mathcal N_{C_K}$ and $D=D_{C_K}$ follow
from Theorem~\ref{thm:functions}. The equivalence of (1) and (2)
is Theorem~\ref{thm:function-positive-splitting}.

The zero-set characterization of $P$-sets
\cite[p.~214]{HenriksenLarsonMartinezWoods1994} gives (2) if and only if (3).
A closed set $K$ is a $P$-set exactly when
\[
 K\subseteq f^{-1}(0)\quad\Longrightarrow\quad
 K\subseteq\operatorname{int}f^{-1}(0)
 \qquad(f\in C(X,\mathbb R)).
\]
In the present notation this condition is precisely $N=D$.

If $N=D$, restriction to $K$ induces an $\ell$-group isomorphism
$C(X,\mathbb R)/D\cong C(K,\mathbb R)$ by the Tietze extension
theorem, so the quotient is Archimedean. Conversely, if
$H=C(X,\mathbb R)/D$ is Archimedean, then
$\operatorname{Inf}_v(H)=\{0\}$. Theorem~\ref{thm:functions}
identifies this ideal with $N/D$, giving $N=D$.
Thus (2) and (4) are equivalent. The existence of a real-linear
section, and hence of an abelian-group section, follows from
Theorem~\ref{thm:function-positive-splitting}.
\end{proof}

For $X=[0,1]$ and $K=\{0\}$, Example~\ref{ex:germs} gives $N\ne D$.
The sequence therefore does not split as a sequence of $\ell$-groups.
The ideals $D$ and $N$ are familiar from continuous-function theory
\cite{GillmanJerison1976}. Related representations by sheaves and
ordered function spaces are developed in \cite{Ball2014,Gehrke2014}.
The role of radicals in local MV-algebras is studied from a different
viewpoint in \cite{CaramelloRusso2017}.

\section{Conclusion}\label{sec:conclusion}

In this paper, we showed that, for every proper lattice filter, the separation axioms
$T_0$, $T_1$, and $T_2$ are equivalent, answering the separation question recorded in \cite{GanLuanDengYang2024}. The zero-neighborhood ideal $N_F$ determines the canonical
Hausdorff reflection, while the polar obstruction ideal $D_F\subseteq N_F$ is the least kernel making
all image radii weak order units. Their quotient laws describe the intermediate
filter-preserving quotients.

For $A=\Gamma(G,u)$, the corresponding polar obstruction ideal $D_C$ is exactly the obstruction
to identifying the filter topology with the subspace topology induced by the ambient
$C$-topology. Thus
\[
\tau_F=\tau_C^G|_{[0,u]}
\iff D_C=\{0\}
\iff \text{every } c\in C \text{ is a weak order unit},
\]
answering Question~2 in \cite{LuanZhaoYang2021}. After quotienting by $D_C$,
the remaining separation defect is the relative infinitesimal ideal. The finite-product
and function-group results make these two obstructions explicit and show that the
residual obstruction also governs the splitting behavior of the associated exact sequence.

These results show that the failure of Hausdorff separation is governed by two successive
algebraic obstructions: a polar obstruction to compatibility with the ambient $C$-topology
and a relative infinitesimal obstruction that remains after the polar quotient.

\section*{Declaration of competing interest}
The authors declare that they have no known competing financial interests or personal relationships that could have appeared to influence the work reported in this paper.

\section*{Acknowledgements}

The first author is supported by  the National Natural Science Foundation of China (12401607, 12571507) and the Basic Research Program of Jiangsu Province (BK20241086). The second author is supported by  the National Natural Science Foundation of China (12471428, 12071033). The third author is supported by  the National Natural Science Foundation of China (12401598) and the Scientific Research Foundation of Nanjing
Institute of Technology (YKJ202445).

\section*{Data availability}
No data were used for the research described in this article.

\end{document}